\documentclass[12pt]{article}
 \newtheorem{theorem}{Theorem}[section]
 \newtheorem{proposition}[theorem]{Proposition}
 \newtheorem{lemma}[theorem]{Lemma}

 \newcommand{\nl}{\newline}

 \newcommand{\dist}{{\rm dist}}

 \newcommand{\R}{{\bf R}}
\newcommand{\C}{{\bf (C)}}

 \newcommand{\diam}{{\rm diam}}

 \newcommand{\diver}{{\rm div}}
 \newcommand{\supp}{{\rm supp}}

 \newcommand{\ia}{({\rm i})}
 \newcommand{\ib}{({\rm ii})}
 \newcommand{\ic}{({\rm iii})}
 \newcommand{\hkp}{{\biggl|\frac{k-p}{p}\biggl|}}

 \title{A unified approach to  improved $L^p$ Hardy inequalities
with best constants}
 \author{
G. Barbatis\footnote{Department of Applied Mathematics,
 University of Crete, 71409 Heraklion, Greece}
 \and S. Filippas\footnotemark[\value{footnote}] \and A. Tertikas
 \footnote{Department of Mathematics,
 University of Crete, 71409 Heraklion, Greece and  \nl
Institute of  Applied and Computational Mathematics,
FORTH, 71110 Heraklion, Greece}
}
 \date{\today}

\begin{document}
 \maketitle

\begin{abstract}
We present a unified approach to  improved $L^p$ Hardy inequalities in $\R^N$.
We consider Hardy potentials that  involve either  the distance from
a point, or  the distance from the boundary, or even  the intermediate case where distance
is taken from a  surface of codimension $1<k<N$. In our main result we add to the right hand
side of the classical Hardy inequality, a weighted $L^p$ norm with optimal weight and
best constant. We also prove  non-homogeneous improved Hardy inequalities, 
where the right hand side involves weighted $L^q$ norms, $q \neq p$.\nl
{\bf AMS Subject Classification: }35J20, 26D10 (46E35, 35P)\nl
{\bf Keywords: } Hardy inequalities, best constants, distance function, weighted norms
\end{abstract}

 \section{Introduction}

The  classical Hardy  inequality  asserts that for  any $p>1$
 \begin{equation}
 \int_{\R^N}|\nabla u|^pdx\geq
 \bigg|\frac{N-p}{p}\bigg|^p\int_{\R^N}\frac{|u|^p}{|x|^p}dx,
 ~~~~~~~u\in C^{\infty}_{c}(\R^N\setminus\{0\}),
 \label{eq:h1}
 \end{equation}
 with  $|\frac{N-p}{p}|^p$ being the best constant,
see for example \cite{HLP}, \cite{OK}, \cite{DH}.
 The best constant remains the same  if $\R^N$ is replaced
  by a domain $\Omega \subset \R^N$ containing the origin.
 Moreover, if $\Omega\subset \R^N$ is a convex domain, possibly unbounded, with
 smooth boundary, and $d(x) =\dist(x,\partial\Omega)$
 the following Hardy inequality:
 \begin{equation}
 \int_{\Omega}|\nabla u|^pdx\geq
 \Bigl( \frac{p-1}{p} \Bigr)^p\int_{\Omega}\frac{|u|^p}{d^p}dx,
 ~~~~~~~u\in C^{\infty}_{c}(\Omega),
 \label{eq:h2}
 \end{equation}
 was recently established, with $(\frac{p-1}{p})^p$ being the best constant, cf
 \cite{MS}, \cite{MMP}. See \cite{OK} for a comprehensive account of Hardy inequalities
and \cite{D} for a review of recent results.

Recently improved versions of (\ref{eq:h1}) and (\ref{eq:h2}) have been obtained.
In  \cite{BV} it is shown  that 
 for  a  bounded domain $\Omega\subset\R^N$ there holds
\begin{equation}
\int_{\Omega}|\nabla u|^2dx-
 \bigg|\frac{N-2}{2}\bigg|^2\int_{\Omega}\frac{u^2}{|x|^2}dx\geq 
c_{\Omega^*}\int_{\Omega}u^2dx,
\label{eq:bv}
\end{equation}
with $c_{\Omega^*}=\Lambda_2(\omega_N/|\Omega|)^{2/N}$, where $\Lambda_2=5.783...$
is the square of the first zero of the Bessel function $J_0$.
It was shown in \cite{FT} that
the optimal constant $c_{\Omega}$ in (\ref{eq:bv}) satisfies 
$c_{\Omega}>c_{\Omega^*}$,
unless $\Omega$ is a ball centered at the origin.
In \cite{GGM} estimate (\ref{eq:bv}) was generalized for $1<p<N$.
It was shown that
when $2\leq p<N$ one can take $c_{\Omega^*}=\Lambda_p(\omega_N/|\Omega|)^{p/N}$
(here $\Lambda_p$ is the first eigenvalue of the $p$-Laplacian in the unit
ball in `$p$-dimensions'); this is not true when $1<p<2$.

In another direction, in  \cite{VZ}, Hilbert space methods were used to
 derive  the following  Improved Hardy-Poincar\'{e}
Inequality
\begin{equation}
\int_{\Omega}|\nabla u|^2dx-
 \bigg(\frac{N-2}{2}\bigg)^2\int_{\Omega}\frac{u^2}{|x|^2}dx
\geq c\left(\int_{\Omega}|\nabla u|^q dx\right)^{2/q}
\label{eq:vz}
\end{equation}
for any $1\leq q<2$.

Analogous results have been obtained in the case of Hardy inequalities
with distance from the boundary. In particular it was proved in
\cite{BM} that for bounded and convex domains there holds
\begin{equation}
\int_{\Omega}|\nabla u|^2dx-
\frac{1}{4}\int_{\Omega}\frac{u^2}{d^2}dx
\geq \frac{1}{4L^2}\int_{\Omega}u^2dx
 \label{eq:bm}
 \end{equation}
and
\begin{equation}
\int_{\Omega}|\nabla u|^2dx-
\frac{1}{4}\int_{\Omega}\frac{u^2}{d^2}dx
\geq \frac{1}{4}\int_{\Omega}\frac{u^2}{d^2(
1 -\log (d/L))^2}dx,
\label{eq:bm1}
\end{equation}
where $L=\diam(\Omega)$.

Hardy inequalities as well as their improved versions have various applications in
the theory of partial differential equations and nonlinear analysis.
They have been useful in the study
of the stability of solutions of semi-linear elliptic and parabolic equations
\cite{PV}, \cite{BV}, \cite{V} as well as in
the existence and asymptotic behavior of the heat equation with
singular potentials, cf \cite{BC}, \cite{CM}, \cite{VZ}; see also \cite{GP} for the
$p$-heat equation. They have also been used to investigate the stability
of eigenvalues in elliptic problems \cite{D,FHT}.

In this work we present a general  approach to
improved Hardy inequalities valid for any $p>1$ and for different choices
of the distance function $d(x)$: besides the two cases above -- distance from a point
and distance from the boundary -- we consider the
more general case where $d(x)$ is the distance of $x\in\Omega$ from a piecewise
smooth surface $K$ of codimension $k$, $1\leq k\leq N$.
 In case $k=N$ we adopt
the convention that $K$ is a point.

In our approach the following geometric assumption on $K$ and $\Omega$ is crucial: if $d(x) = \dist(x,K)$ then the following inequality should
hold in the weak sense:
\[p\neq k, \qquad\qquad\Delta_p d^{\frac{p-k}{p-1}}\leq 0, 
\quad\quad\mbox{ in $\Omega\setminus K$.}\qquad\hfill\hspace{2cm}  \C
\]
Here $\Delta_p$ denotes the usual $p$-Laplace operator,
$\Delta_p w=\diver(|\nabla w|^{p-2}\nabla w)$.
This condition is analyzed in detail in Section 2. Here we simply 
 note that (C) is always satisfied when $k=N$ and $d(x)$ measures
the distance from a point as well as when
$k=1$, $\Omega$ is convex and $d(x)$ is
the distance from $K=\partial\Omega$.
Condition (C) can be interpreted as a higher-codimension analogue
of the usual convexity condition that appears in Hardy's inequality
when $k=1$ and
$K=\partial\Omega$; cf (\ref{eq:h2}).

In order to describe our  results we introduce the function
 \[X(t) = -1/ \log t,\quad t\in (0,1).\]
Our main theorem then is the following:

{\bf Theorem A (Improved Hardy Inequality)}
{\em Let $\Omega$ be a domain in $\R^N$ and $K$ 
a piecewise smooth surface of codimension $k$, $k=1,\ldots,N$.
Suppose that  $\sup_{x\in\Omega}d(x,K) < \infty$
and       condition (C) is satisfied. Then \nl
(1) There exists a positive constant
$D_0=D_0(k,p)\geq \sup_{x\in\Omega}d(x,K)$ such that for any $D \geq  D_0$
and all  $u \in W^{1,p}_0(\Omega\setminus K)$ there holds
 \begin{equation}
 \int_{\Omega}|\nabla u|^pdx -\hkp^p\int_{\Omega}\frac{|u|^p}{d^p}dx
 \geq  \frac{p-1}{2p}\hkp^{p-2}
\int_{\Omega}\frac{|u|^p}{d^p}X^2(d/D)dx.
 \label{improved}
 \end{equation}
If in addition $2 \leq p <k$, then we can take
$D_0 = \sup_{x \in \Omega} d(x,K)$.\nl
(2) Both constants appearing in (\ref{improved}) as well as 
the exponent two in  $X^2$ are optimal
in either of the following cases:
\begin{eqnarray*}
(a) && \mbox{ $k=N$ and
$K=\{0\} \subset \Omega$;}\\
(b) &&\mbox{ $k=1$ and $K =\partial \Omega$.} \\
(c) &&\mbox{ $2\leq k\leq N-1$ and $\Omega\cap K\neq\emptyset$;}
\end{eqnarray*}}

The optimality of the constants and the exponent
is meant in the following sense:
\[
\biggl|\frac{p-k}{p}\biggr|^p=\inf\Bigl\{\int_{\Omega}|\nabla u|^pdx,\;\;
 \int_{\Omega}\frac{|u|^p}{d^p}dx =1\Bigr\}.
\]
Further, if $\gamma<2$, then, no matter how large $D$ is,
there is no $c>0$ such that
\[
\int_{\Omega}|\nabla u|^pdx - \biggl|\frac{p-k}{p}\biggr|^p
\int_{\Omega}\frac{|u|^p}{d^p}dx
 \geq  c\int_{\Omega}\frac{|u|^p}{d^p}X^{\gamma}(d/D)dx;
\]
and finally, for any $D\geq D_0$,
\[
\frac{p-1}{2p}\biggl|\frac{p-k}{p}\biggr|^{p-2}=\inf\Bigl\{
\int_{\Omega}|\nabla u|^pdx -  \hkp^p\int_{\Omega}\frac{|u|^p}{d^p}dx
,\;   \int_{\Omega}\frac{|u|^p}{d^p}X^2(d/D)dx=1\Bigr\}.
\]

A few remarks are in order:

1. The assumption $D \geq D_0$ is only
necessary in order to obtain the precise constant $(p-1)/(2p)|(k-p)/p|^{p-2}$.
We can take any $D>\sup_{x\in\Omega}d(x,K)$ at the expense of having
a smaller constant  $c=c(p,k,D)$ in the right hand side of (\ref{improved}).

2. The logarithmic correction in the right hand side is independent of $p>1$.
Also it is worth pointing out that the constant of the Improved Hardy Inequality
depends only on $p$ and $k$ and not on $K$, the dimension $N$, or $\Omega$.
This is in contrast to the Improved Hardy
Inequalities which involve the unweighted $L^p$ norm in the right
hand side (see e.g. (\ref{eq:bv}), (\ref{eq:bm})).

3. A simple density argument shows that if $p<k$ then
$W^{1,p}_0(\Omega\setminus K)=W^{1,p}_0(\Omega)$.

4. We only assume that $\dist(x,K)$ is bounded on $\Omega$, not
that $\Omega$ itself is bounded.

In case $p=2$ and $k=1$ or $N$, part (1) of Theorem A has been
obtained in \cite{BM,BV} by a different method. We are aware of
 very few results in the literature for  $1<k<N$, concerning 
even the simple Hardy inequality with best constant; for the case $p=2$
 see \cite{D, DM}, and
 \cite{M} Section 2.1.6.

We present two different approaches to the Improved Hardy Inequality.
The first  is based on a suitable change of variables \cite{BM,BV,GGM,M}. While
this method does not yield the optimal constant in the right hand side
 of (\ref{improved}),
it has the advantage that it easily leads to non-homogeneous improved
Hardy inequalities. We note that in this method the arguments
used for $1<p<2$ differ from those used for $p\geq 2$.
The second approach is based on the careful choice of a suitable
vector field and an elementary integral inequality and is the one that
gives the sharp constants.
It is remarkable that condition (C) comes up naturally in both approaches.

It is well known that for $k=N$ (distance from a point) there is no
Hardy inequality if  $p=N$. More generally there is no Hardy inequality 
if $p=k$, $1 \leq kj \leq N$. For  that case  we provide a substitute for Hardy inequality
with optimal weight and best constant; see Theorems \ref{thm:os} and  \ref{thm:opos}.

We next consider non-homogeneous improved Hardy inequalities which
involve $L^q$ norms with  $q \neq p$, in the right hand side.
In this direction we have the following:

 {\bf Theorem B (Improved Hardy-Poincar\'{e} Inequality)}
{\em Let $\Omega$ be a domain in $\R^N$ and $K$ 
a piecewise smooth surface of codimension $k$, $k=1,\ldots,N$.
Suppose that  $\sup_{x\in\Omega}\dist (x,K) < \infty$
and       condition (C) is satisfied. Then \nl
(1) For any $D>\sup_{x\in\Omega}\dist (x,K)$, $1 \leq q <p$ and $\beta>1+q/p$ 
there exists a positive constant $c>0$ such that for all $ u\in  W^{1,p}_0(\Omega\setminus K)$
there holds:
 \begin{equation}
 \int_{\Omega}|\nabla u|^pdx - \hkp^p\int_{\Omega}\frac{|u|^p}{d^p}dx
 \geq
 c\left(\int_{\Omega}|\nabla u|^q d^{k(-1+q/p)}X^{\beta}(d/D)dx\right)^{p/q}.
 \label{eq:hp}
 \end{equation}

(2)The estimate is sharp in the sense that 
 the exponent of $X$ in the right hand side of (\ref{eq:hp})
cannot be smaller than $1+q/p$, in either of the cases (a), (b), (c) of
part (2) of Theorem A.}

For $p=2$ and $k=N$ this strengthens inequality (\ref{eq:vz}).

We next consider Improved Hardy Sobolev inequalities. 
 Let $K=\{x\in\R^N\,|\,x_1=x_2=\ldots =x_k=0\}$, $1\leq k\leq N-1$.
Then, in \cite{M} the following inequality is established (see Corollary 3,
section 2.1.6) for any $2<q \leq \frac{2N}{N-2}$:
\begin{equation}
 \int_{\Omega}|\nabla u|^2dx - \biggl| \frac{k-2}{2}\biggr|^2\int_{\Omega}\frac{u^2}{d^2}dx
 \\
 \geq c\left(\int_{\Omega}|u|^q d^{-q-N+Nq/2} dx
 \right)^{2/q},
 \label{eq:maz1}
 \end{equation}
for any $u \in C_c^{\infty}(\R^N \setminus K)$. The question was posed
in \cite{M} whether an analogue result holds for $p \neq 2$.

For $k=N$ that is $K=\{0\} \in \Omega$, a bounded domain in $ \R^N$ an analogous inequality is shown
 in \cite{BV}, valid  any $2 \leq q <\frac{2N}{N-2}$:
\[
 \int_{\Omega}|\nabla u|^2dx - \biggl|\frac{N-2}{2}\biggr|^2\int_{\Omega}\frac{u^2}{|x|^2}dx
\geq c \left( \int_{\Omega} |u|^q dx \right)^{2/q}.
\]

Our result reads:

{\bf Theorem C (Improved Hardy-Sobolev Inequality)}

{  \em (1)  Let $K=\{x\in\R^N\,|\,x_1=x_2=\ldots =x_k=0\}$, $1\leq k\leq N-1$. Assume that
 $2 \leq p<N$ and $p < q \leq Np/(N-p)$.
  Then there exists a constant $c>0$ such that for all $u \in W_0^{1,p}(\R^N \setminus K)$
there holds
\begin{equation}
 \int_{\R^N}|\nabla u|^pdx -  \hkp^p  \int_{\R^N}\frac{|u|^p}{d^p}dx
 \\
 \geq c\left(\int_{\R^N}|u|^q d^{-q-N+Nq/p} dx
 \right)^{p/q}.
 \label{eq:hs1}
 \end{equation}
(2) Let  $k=N$, that is $K=\{0\}\in \Omega$, a bounded domain in $\R^N$. Assume that $1<p<N$ and
$p \leq  q < Np/(N-p)$.
 Then there exists a constant $c>0$ such that for all  $u \in W_0^{1,p}(\Omega)$
\begin{equation}
 \int_{\Omega}|\nabla u|^pdx - \hkp^p\int_{\Omega}\frac{|u|^p}{d^p}dx
 \\
 \geq c\left(\int_{\Omega}|u|^q d^{-q-N+Nq/p}X^{1+q/p}(d/D)dx
 \right)^{p/q}.
 \label{eq:hs2}
 \end{equation}
Inequality (\ref{eq:hs2}) is optimal in the sense that
 $X^{1+q/p}$ cannot be replaced by a smaller power of $X$.}

A simple scaling argument shows that the exponent of $d$ in (\ref{eq:hs1}) is optimal.
Hence it comes as a remarkable fact that the case $k=N$ is different from the case $k<N$.
It is an open question whether (\ref{eq:hs2}) remains true in the critical
case $q=Np/(N-p)$. One can see that for $q=Np/(N-p)$ one cannot have
an inequality (\ref{eq:hs2}) without the presence of the logarithmic
correction. In fact one cannot even  have the weak $L^{Np/(N-p)}$ norm in the
right hand side; see
Proposition \ref{prop:weak}. On the other hand inequality (\ref{eq:hs2})
is true in the critical case if we replace $X^{1+q/p}$ by $X^{2q/p}$. This
last result is contained in  Theorem \ref{thm:sobboun} where an  inequality 
weaker than (\ref{eq:hs1}) and (\ref{eq:hs2}) is
shown valid  for   $k \leq N$ and non-affine $K$.

 The structure of the paper is as follows: in Section 2 we
discuss the geometric assumptions on $\Omega$ and $K$; in particular we
provide specific examples for which condition (C) is satisfied.  Section 3
contains our first approach to Improved hardy inequalities, whereas Section 4
is devoted to the vector field approach which yields the best constants.
In Section 5 we prove the optimality of the constants involved in Theorem A.
Finally in Section 6, we use the results of Section 3 to obtain non-homogeneous
inequalities.

{\bf Acknowledgment.} We would like to thank E.M. Harrell II for
a useful discussion and for bringing to our attention  inequality (\ref{eq:harrell})
below.





 \setcounter{equation}{0}
 \section{The geometry of $K$ and $\Omega$}

 In this section we shall introduce the main geometric assumptions
 concerning $K$ and $\Omega$, and we will  fix
 some notational conventions.
 Throughout this work $\Omega$ is a domain in $\R^N$ and
 $K$ is a piecewise smooth closed and connected surface  of codimension
 $k=2,3...N-1$. We also allow for the two extreme cases $k=1$ or $N$, with the following
 convention: If $k=N$ then $K$ is reduced to a point, say the origin.
 If $k=1$, then we take  $K$ to be  the
 boundary of $\Omega$, that is $K= \partial \Omega$.

 In all cases we define  the distance function $d(x)$ by
 \[
 d(x) =  \dist (x,K),\qquad x \in \Omega.
 \]
 Hence for $k=N$ we have $d(x) = |x|$, whereas for  $k=1$, $d(x)$ is
 the distance from the boundary of $\Omega$.
 Let us note that $d(x)$ is 
 a Lipschitz continuous function with $|\nabla d|=1$ a.e..

We now come to our main geometric assumption on $K$ and $\Omega$,
expressed in terms of the distance function $d$.
We introduce the following geometric condition: 
\[p\neq k\;\;\mbox{ and }\;\;\Delta_p d^{\frac{p-k}{p-1}}\leq 0 \;\;
\mbox{ on $\Omega\setminus K$}.\qquad\hfill\hspace{2cm} \C
\]
Formal calculations give
\[
\Delta_p d^{\frac{p-k}{p-1}}=\frac{p-k}{p-1}\biggl|\frac{p-k}{p-1}
\biggr|^{p-2}d^{-k}(d\Delta d+(1-k)|\nabla d|^2),\\
\]
so that, since $|\nabla d|=1$ a.e., an equivalent formulation of (C) is
\[(p-k)(d\Delta d+1-k)\leq 0 \;
\mbox{ on $\;\Omega\setminus K$}.\]

The precise meaning of the above condition is the following: we consider the linear functional
 \begin{eqnarray*}
A[\phi]&:=& -\int_{\Omega}  |\nabla d|^2 \phi dx
 -\int_{\Omega} d \, \nabla d \cdot  \nabla \phi  \,  dx +(1-k)\int_{\Omega} \phi dx \\
&=&-\int_{\Omega} d \, \nabla d \cdot  \nabla \phi  \,  dx -k\int_{\Omega} \phi dx,\quad
\phi\in C^1_c(\Omega\setminus K), 
\end{eqnarray*}
and require that for all non-negative $\phi\in C^1_c(\Omega\setminus K)$ there holds
$(p-k)A[\phi]\leq 0$. 
In this context, and in order to simplify our notation, we shall
use the expression
\[\int_{\Omega}(d\Delta d+1-k)\phi \,dx,\quad \phi\in C^1_c(\Omega\setminus K),\]
to denote the functional $A[\phi]$. This allows us to perform formal
integrations by parts as if $\Delta d$ were a locally integrable function in $\Omega$.
Taking for instance $\phi=\psi/d$ in the definition above we obtain the relation
\[\int_{\Omega}\psi\Delta d\,dx=-\int_{\Omega}\nabla\psi\cdot\nabla d\,dx,\qquad
\psi\in C^1_c(\Omega\setminus K).\]
This also justifies the following convention: assuming that (C) is satisfied, we define:
\[\int_{\Omega}\left|d\Delta d+1-k\right| \phi \,dx=\left\{
\begin{array}{rl} 
A[\phi], & \mbox{if $p<k$} \\
-A[\phi], & \mbox{if $p>k$;}
\end{array}\right.\]
this is
a positive functional on $C^1_c(\Omega\setminus K)$ and it is then easily seen that
\[\left|\int_{\Omega}\left(d\Delta d+1-k\right) \phi \,dx\right|
\leq \int_{\Omega}\left|d\Delta d+1-k\right| |\phi| \,dx,\quad
\phi\in C^1_c(\Omega\setminus K).\]

We next present some examples in which condition (C) is satisfied.
The first two concern the cases $k=1$ and $k=N$, which are the most popular in the literature.
Then we consider the intermediate cases $2\leq k\leq N-1$. One then is lead to
rather special  assumptions on $K$ and $\Omega$. This is not due to lack of pairs
$(K,\Omega)$ that satisfy (C); indeed, it is easy to see that
given any $K$ one can always find an $\Omega$ such that (C) is satisfied: simply take
$\Omega$ to be any domain contained in the set
\[\{x\in\R^N \, :\, d\Delta d+1-k \geq 0\: \mbox{ (or $\leq 0$)}\}.\]
An analytical description of such sets $\Omega$ is possible only after extra
assumptions on $K$.

{\bf Example 1. }Let  $k=N$ so that $K=\{0\}$.
Then $d(x) = |x|$ and $\Delta d^{2-N}=0$ away from $x=0$,
hence condition (C) is satisfied for any $1<p<\infty$ and any $\Omega\subset\R^N$.

{\bf Example 2. }Suppose that $k=1$, so that $K= \partial \Omega$. Then (C) is satisfied
for all $1<p<\infty$ provided we make the additional assumption that $\Omega$ is {\em convex}. 
To see this we first claim that $d(x)$, $x \in \Omega$, is a concave function. Indeed,
 let $0< \lambda <1$, and  $x$, $y$, $z= \lambda x + (1-\lambda) y$ be three points
  contained
 in $\Omega$. Let $z_0 \in \partial \Omega$ be a point that realizes the distance
 for  $z$, that is, $d(z) = |z-z_0|$. We denote by $T_{z_0}$ the hyperplane
 that contains $z_0$ and is orthogonal to the vector $z-z_0$.
 We also let $x_0$ and $y_0$ be the orthogonal projections
 of $x$ and $y$ onto $T_{z_0}$ respectively. It then follows by the convexity
 of $\Omega$ and a simple similarity argument that
 \[
 d(z) = |z-z_0| = \lambda |x-x_0| + (1- \lambda) |y - y_0|
 \geq  \lambda d(x) + (1- \lambda) d(y),
 \]
and the claim is proved. Since $d(x)$ is concave we conclude from
Theorem 6.3.2 of \cite{EG} that $\Delta d$ is non-positive
in the weak sense; more precisely there exists a non-negative Radon
measure $d\mu$ on $\Omega$ satisfying
\begin{equation}
\int_{\Omega}\nabla\psi\cdot\nabla d\, dx=\int_{\Omega}\psi d\mu,\quad\psi\in C^1_c(\Omega);
\label{eq:ibp}
\end{equation}
In particular, taking as test function $\psi =\phi d$, we see that
$A[\phi]\leq 0$, that is condition (C) is satisfied.

Let us now consider the intermediate  cases  $2 \leq k \leq N-1$.

{\bf Example 3. }If $K$ is affine, $K\equiv\R^{N-k}$, then condition (C)
is satisfied for all $1<p<\infty$ without any restriction on $\Omega$.\nl
Indeed, changing coordinates if necessary, we see by a
direct computation that
\[
\Delta_p d^{\frac{p-k}{p-1}} = 0,\quad x\in\R^N\setminus K.
\]
Further, if $p>k$ and $K$ is the union of affine sets,
\[K=\cup_{i\in I}K_i,\]
then (C) is also satisfied, again with no restriction on $\Omega$.
To see this consider the functions $d_i(x)=\dist(x,K_i)$, $i\in I$.
We have seen that $d_i^{(p-k)/(p-1)}$ is $p$-harmonic. But
\[ d^{\frac{p-k}{p-1}}(x)=\inf_{i\in I}d_i^{\frac{p-k}{p-1}}(x),\quad
x\in\Omega\setminus K\]
and hence $d^{\frac{p-k}{p-1}}$ is $p$-super-harmonic by
the comparison principle for the $p$-Laplacian, see \cite{HKM}.
Alternatively, observing that
\[
\Delta_p d^{\frac{p-k}{p-1}}= \frac{p-k}{p-1} \, \biggl|\frac{p-k}{p-1}
\biggr|^{p-2}   \frac{1}{2-k}  \Delta d^{2-k},
\]
we may use the corresponding  principle for the Laplacian. (When $k=2$
we replace $\frac{1}{2-k} \Delta d^{2-k}$ by $\Delta \log d$.)

{\em Definition} Let $E \subset \R^N$  be an affine set of codimension $k-1$ and
$V\subset E$ be a convex domain (i.e. connected and
open in the topology of $E$) and let $K = \partial_E V$.
The cylinder $V\times E^{\perp}$ is called the {\em inner canal} of $K$;
the cylinder $(E\setminus \overline{V})\times E^{\perp}$ is called the
{\em outer canal} of $K$. (See also \cite{S})

{\bf Example 4. }(i) If $p>k$ and $\Omega$ is contained in the inner canal of $K$
then (C) is satisfied; (ii) If $p<k$ and $\Omega$ is contained in the outer canal of $K$
then (C) is satisfied.\nl
To see (i) let $\{T_y \, |\, y\in\partial_E V\}$ be the family of hyperplanes in $E$
which are tangent to $K$ (if $K$ is not smooth we take the supporting hyperplanes instead).
If $\Omega$ is contained in the inner canal of $\Omega$ then
\[d(x)=\inf_{y\in\partial_E V}\dist(x,T_y),\quad x\in\Omega\]
and we are back in the situation of Example 3. \nl
To prove (ii) we use a different argument.
We write any $x \in \R^N$ as $x = (y,z)$ with
 $y \in E\equiv \R^{N-k+1}$ and $z \in \R^{k-1}$; that is, the projection
 of $x$ onto $E$ is the point $(y,0)$. We then have
 \begin{equation}
 d^2(x) = \tilde{d}^2(y) + |z|^2, \qquad\tilde{d}(y) = \dist((y,0), K).
\label{eq:cal} 
\end{equation}
 Differentiating twice with respect to $z_i$, $i=1,2,\ldots,k-1$,
 and summing up over $i$ we obtain
\begin{equation}
|\nabla_{z} d|^2 + d \Delta_z d \, = k-1.
\label{eq:dz}
\end{equation}
Differentiating (\ref{eq:cal})
with respect to $y_i$, $i=1,2,\dots,N-k+1$,
we obtain in a similar way
 \begin{equation}
 |\nabla_{y} d|^2 + d \Delta_y d \, =
 |\nabla_y \tilde{d}|^2 +  \tilde{d} \Delta_y  \tilde{d} =
 1 +    \tilde{d} \Delta_y  \tilde{d}.
 \label{eq:dy}
 \end{equation}
Adding (\ref{eq:dz}) and (\ref{eq:dy}) we conclude that
 \[
  d\Delta d+(1-k)|\nabla d|^2 =  \tilde{d} \Delta_y  \tilde{d}.
 \]
Since $\tilde{d}$ is the distance function in  $E \equiv \R^{N-k+1}$
and $V \subset E$ is a convex domain, we have, as in Example 2,
that $ \Delta_y  \tilde{d} \geq 0$ if $y \in V^c$. Hence (C)
is satisfied in this case.\nl
We point out that if a domain $\Omega$ satisfies $\overline{\Omega}
\cap K\neq\emptyset$ (so that $d^{-1}$ is singular in  $\Omega$)
then for it to be contained in either the
inner or the outer canal of $K$ it is necessary that 
$K \cap \partial\Omega \neq \emptyset$.

Our fifth example combines ideas from the last two ones.

{\bf Example 5. }Assume that $p>k$ and that $\Omega$ is contained in the
inner  canal of $L= \partial V$. Let $K$ be a polytope contained in $V$ and
having its vertices on $L$. Then condition (C) is satisfied. To see this let
$F_i$, $i=1,\ldots,L$, be the faces of $K$. Our assumption on $K$ and $\Omega$
imply that the distance of any $x\in\Omega$ from a face $F_i$ is realized at a point $y\in F_i$
which is on the interior of $F_i$, that is, the distance is not realized at
vertices, edges etc. Hence 
\[\Delta_p d^{\frac{p-k}{p-1}}=0,\quad x\in\Omega\setminus F_i,\quad d_i(x)=\dist(x,F_i),\]
and the comparison argument of Example 3 goes through.


 \setcounter{equation}{0}
 \section{The Improved Hardy inequality}

 In this section we give a first proof of the improved Hardy
inequality and also obtain some inequalities which
will be of use in Section 6. We start
 with some elementary pointwise inequalities.
 \begin{lemma}
 For any $1<p<\infty$ there exists a constant
 $c=c(p)>0$ such that for all $a,b\in\R^N$ we have:
 \begin{eqnarray*}
\ia &&\mbox{ if $1<p<2$ then}\\
 && {|a-b|^p-|a|^p\geq c\frac{|b|^2}{(|a|+|b|)^{2-p}}
 -p|a|^{p-2}a\cdot b.} \\
 \ib &&~ \mbox{if $p\geq 2$ then}\\
 &&{\rm (a)}\hspace{.8cm}
 {|a-b|^p-|a|^p\geq  c|a|^{p-2}|b|^2  -p|a|^{p-2}a\cdot b\, ;}\\
 &&{\rm (b)}\hspace{.8cm}
 {|a-b|^p-|a|^p\geq  c|b|^p  -p|a|^{p-2}a\cdot b\, ;}
 \end{eqnarray*}
 \label{lem:lin}
 \end{lemma}
 {\em Proof. }Parts $\ia$ and $\ib$(b) are contained in
 Lemma 4.2 of \cite{L}. Hence, we only  prove $\ib$(a).

If $|b| \geq \frac{1}{2} |a|$ the inequality follows from   $\ib$(b).
Suppose now that  $|b| < \frac{1}{2} |a|$; then $|a -\xi b| \geq \frac{1}{2}|a|$
for all $\xi\in(0,1)$. Hence, taking the Taylor expansion
of $f(t) = |a-bt|^p$ around $t=0$ we have,
\begin{eqnarray*}
|a-b|^p & =&  |a|^p-p|a|^{p-2}a\cdot b + \frac{p(p-2)}{2}
|a -\xi b|^{p-4}((a-\xi b)\cdot b)^2 +\\
&&+p|a-\xi b|^{p-2}|b|^2 \hspace{5.5cm}\mbox{ (some $\xi\in (0,1)$)}\\
&\geq&  |a|^p-p|a|^{p-2}a\cdot b + \frac{p}{2^{p-2}} |a|^{p-2}|b|^2.
\end{eqnarray*}
$\hfill //$

We next prove an auxiliary inequality that will be used in the sequel. 
Let us first recall that 
\[X(s)=-\frac{1}{\log s},   \qquad s \in(0,1).\]
Note that if $D > \sup_{x \in \Omega} d(x)$ then 
\[0\leq  X(d(x)/D)\leq M,\qquad x\in\Omega,\]
for a suitable positive constant $M=M(D)$.
Furthermore, we shall often use the relation
\begin{equation}
\frac{d}{dr}X^{\beta}=\beta\frac{X^{\beta+1}}{r},
\label{eq:de}
\end{equation}
as well as its integral version
\begin{equation}
\int_{s_1}^{s_2}r^{-1}X^{\beta +1}(r)dr=\frac{1}{\beta}
[X^{\beta}(s_2)-X^{\beta}(s_1)].
\label{eq:int}
\end{equation}

We next prove the following

 \begin{lemma}
 Let $D> \sup_{x \in \Omega} d(x,K)$ and
  $\alpha\neq 1$. Then
 \[
 \left( \frac{|\alpha-1|}{p}\right)^p\int_{\Omega}|v|^pd^{-k}X^{\alpha}(d/D)dx\leq
 \int_{\Omega}|\nabla v|^pd^{p-k}X^{\alpha-p}(d/D)dx
 \]
 \begin{equation}
  +\left(\frac{|\alpha-1|}{p}\right)^{p-1}\int_{\Omega}|v|^p  d^{-k}
 \left|d \Delta d+1-k\right| X^{\alpha-1}(d/D) dx
 \label{eq:sh1}
 \end{equation}
 for all $v\in C^{\infty}_c(\Omega \setminus K)$.
 \label{lem:tra}
 \end{lemma}

 {\em Proof. }We prove (\ref{eq:sh1}) for $D=1$, the general case following
 by scaling. Recalling (\ref{eq:de}) we have
 \begin{eqnarray*}
 &&\hspace{-1cm}\int_{\Omega}|v|^pd^{-k}X^{\alpha}(d)dx=\\
&=&\frac{1}{\alpha-1}\int_{\Omega}|v|^p d^{1-k} \nabla d\cdot\nabla X^{\alpha-1}(d)dx\\
 &=&-\frac{p}{\alpha-1}\int_{\Omega}|v|^{p-2}vX^{\alpha-1}(d) d^{1-k}
 \nabla v\cdot\nabla d\: dx \\
 & &-\frac{1}{\alpha -1}\int_{\Omega}|v|^p  d^{-k} (d \Delta d + (1-k) |\nabla d|^2) 
X^{\alpha-1}(d)\,dx\\
 &\leq&\frac{p}{|\alpha-1|}\left(\int_{\Omega}|\nabla v|^pd^{p-k}X^{\alpha-p}(d)dx\right)^{1/p}
 \left(\int_{\Omega}|v|^pd^{-k}X^{\alpha}(d)dx\right)^{(p-1)/p}+\\
 &&+\frac{1}{|\alpha -1|}\int_{\Omega}|v|^p
   d^{-k}\left| d \Delta d +1-k\right|
   X^{\alpha-1}(d) \, dx.
 \end{eqnarray*}
 Hence we have an estimate of the form
 \[B\leq \theta B^{(p-1)/p}\Gamma^{1/p}+A, ~~~~~~~~~\theta = \frac{p}{|\alpha -1|}.\]
 Combining this with the relation
 \[B^{(p-1)/p}\Gamma^{1/p}\leq\frac{\epsilon(p-1)}{p}B+\frac{\epsilon^{-(p-1)}}{p}\Gamma,\]
 and taking $\epsilon=\theta^{-1}$ we obtain
 $\theta^{-p}B\leq\Gamma +p\theta^{-p}A$, which is the required
 inequality. $\hfill$//

Throughout the paper we will use the notation
\begin{equation}
I[u]=\int_{\Omega}|\nabla u|^pdx-\Bigl|\frac{p-k}{p}\Bigr|^p\int_{\Omega}
\frac{|u|^p}{d^p}dx,\quad u\in W^{1,p}_0(\Omega\setminus K).
\label{eq:iu}
\end{equation}

Our starting point is the following lower estimate on $I[u]$:
 \begin{lemma}
 Let $u\in W^{1,p}_0(\Omega\setminus K)$ be given
 and set $v(x)=u(x)d^{H}(x)$, $H=(k-p)/p$.
 There exists a constant $c=c(p)>0$ such that:
 \begin{eqnarray}
 \ia && \mbox{if $1< p<2$ then }\nonumber\\
 &&\hspace{-.6cm} I[u]\geq c\int_{\Omega}\frac
 {|\nabla v|^2d^{2-k}}{\left(|Hv|+|d\nabla v|\right)^{2-p}}dx  \nonumber \\
 &&\hspace{0.5cm} +H|H|^{p-2}\int_{\Omega}|v|^p  d^{-k}
 (d\Delta d+1-k) dx;\label{eq:a81} \\
 \ib && \mbox{if $2\leq p<\infty$ then }\nonumber\\
 &&({\rm a})\hspace{.8cm} I[u] \geq \,  c |H|^{p-2}
 \int_{\Omega}|\nabla v|^2|v|^{p-2}
 d^{2-k} \,dx  \nonumber  \\
 && \hspace{2.6cm} +H|H|^{p-2}\int_{\Omega}|v|^p d^{-k}
(d\Delta d+1-k)  dx\label{eq:a71};\\
 &&({\rm b})\hspace{.8cm} I[u] \geq \,  c\int_{\Omega}|\nabla v|^pd^{p-k}dx  \nonumber \\
 && \hspace{2.6cm} +H|H|^{p-2}\int_{\Omega}|v|^p     d^{-k}
 (d\Delta d+1-k)  dx;\label{eq:aa11}
 \end{eqnarray}
 \label{lem:fun1}
 \end{lemma}
 {\em Proof. }It is straight forward to see that
 \[ |\nabla u|^p- |H|^p\frac{|u|^p}{d^p}=d^{-k}
 \left(\left|Hv\nabla d-d\nabla v
 \right|^p-|Hv|^p\right);\]
 to estimate the right hand side  we use the corresponding inequalities of
 Lemma \ref{lem:lin} with $a=Hv\nabla d$ and $b=d\nabla v$. The expression
 $-p\int_{\Omega}d^{-k}|a|^{p-2}a\cdot b$ appears in all three cases
 and is equal to $H|H|^{p-2}\int_{\Omega}|v|^p d^{-k}
 (d\Delta d+1-k)  dx$ as can be seen  by an integration by parts. The
 stated estimates then follow at once. $\hfill //$

It should be noted that if condition (C) is satisfied then the 
common term that appears in the right hand side of the three inequalities
of the last lemma is equal to 
\[|H|^{p-1}\int_{\Omega}|v|^p d^{-k}
\left|d\Delta d+1-k\right| dx,\]
and, in particular, is non-negative.

We next prove the improved Hardy inequality for $1<p<2$.
\begin{proposition}
Let $1<p<2$. Given $u\in W^{1,p}_0(\Omega \setminus K)$ we set
$v(x)=u(x)d^{H}(x)$, $H=(k-p)/p$. If condition (C) is satisfied
then there exist  constants $c_i=c_i(p,k,D)>0$, $i=1,2$ such that
 \begin{eqnarray}
 I[u]\!  &\! \geq& c_1 \left(\int_{\Omega} |\nabla v|^p d^{p-k}
 X^{2-p}(d/D) \, dx+  \int_{\Omega}  | v|^p d^{-k}
 \left|d\Delta d+1-k\right| \, dx\right)
 \nonumber  \\
 &\geq &\! c_2  \int_{\Omega}
 \frac{|u|^p}{d^p}X^2(d/D) \, dx.
 \label{eq:1p2}
 \end{eqnarray}
 \label{prop:1p2}
 \end{proposition}
 {\em Proof:}
 We may assume that $D=1$, the general case following by scaling.
 To simplify the  subsequent calculations we set
 \[A_1=\int_{\Omega}\frac
 {|\nabla v|^2d^{2-k}}{\left(|Hv|+|d\nabla v|
 \right)^{2-p}}dx,\quad A_2=\int_{\Omega}|v|^pd^{-k}X^2(d/D)dx\]
 \[ A_3=\int_{\Omega}|\nabla v|^pd^{p-k}X^{2-p}(d/D)dx,\quad
 A_4=\int_{\Omega}|v|^p d^{-k} |d\Delta d+1-k|dx.\]
 Note the all $A_i$'s are positive and homogeneous
 of degree $p$ in $v$.  H\"{o}lder's inequality and  elementary estimates yield
 \begin{eqnarray*}
 A_3 & =& \int_{\Omega} \frac{|\nabla v|^2d^{p(2-k)/2}}{\left(|Hv|+|d\nabla v|
 \right)^{p(2-p)/2}} \,\cdot \, \left(|Hv|+|d\nabla v|
 \right)^{p(2-p)/2} d^{-k(2-p)/2}  X^{2-p} dx \\
 &\leq& A_1^{p/2}\left(\int_{\Omega}\left(|Hv|+
 |d\nabla v|\right)^pd^{-k}X^2 dx\right)^{(2-p)/2}\\
 &\leq&c(p,k) \,  A_1^{p/2}\left(\int_{\Omega}|v|^p d^{-k}X^2dx
 +\int_{\Omega}|\nabla v|^p d^{p-k}X^2 dx\right)^{(2-p)/2}\\
 &\leq&c(p,k) \,  A_1^{p/2}(A_2+A_3)^{(2-p)/2},
 \end{eqnarray*}
 that is
 \begin{equation}
 A_1\geq c(p,k) \frac{A_3^{2/p}}{(A_2+A_3)^{(2-p)/p}}.
 \label{eq:s1}
 \end{equation}
 It  follows from Lemma \ref{lem:fun1}(i) that
 \begin{equation}
 I[u] \geq c(k,p)(A_1 + A_4).
 \label{eq:s2}
 \end{equation}
We also have from Lemma \ref{lem:tra} (with $\alpha = 2$),
 \begin{equation}
 A_2\leq c(p,k) (A_3+A_4).
 \label{eq:s3}
 \end{equation}
 Combining (\ref{eq:s1}), (\ref{eq:s2}) and (\ref{eq:s3})
 we obtain
 \begin{eqnarray*}
 I[u]&\geq&c\left(\frac{A_3^{2/p}}{(A_2+A_3)^{(2-p)/p}}
 +A_4\right)\nonumber\\
 &\geq&c\left(\frac{A_3^{2/p}}{(A_3+A_4)^{(2-p)/p}}
 +A_4\right)\nonumber\\
 &\geq&c(A_3+A_4),
 \end{eqnarray*}
 which is the first inequality in (\ref{eq:1p2}). Using once more (\ref{eq:s3}) we
 have
 \[
 I[u] \geq c(A_3 + A_4) \geq cA_2 = c
  \int_{\Omega}
 \frac{|u|^p}{d^p}X^2 \, dx,
 \]
 and the proof of (\ref{eq:1p2}) is complete.$\hfill //$

We now consider the complementary case $p\geq 2$.
\begin{proposition}
Let $p \geq 2$. Given $u\in W^{1,p}_0(\Omega \setminus K)$
we set $v(x)=u(x)d^{H}$, $H=(k-p)/p$. If condition (C)
is satisfied then there exists a constant $c=c(p,k,D)>0$,  such that
 \begin{equation}
 I[u]
 \geq  c  \int_{\Omega}
 \frac{|u|^p}{d^p}X^2(d/D) \, dx.
 \label{eq:p2}
 \end{equation}
 \label{prop:p2}
 \end{proposition}
{\em Proof.} We will use the additional
change of variables $w = |v|^{p/2}$. 
It follows from Lemma \ref{lem:fun1}(iia) that
 \begin{eqnarray*}
  I[u] & \geq &  c
 \int_{\Omega}|\nabla v|^2|v|^{p-2}
 d^{2-k} \,dx +
 c \int_{\Omega}|v|^p d^{-k}
 \left|d\Delta d+1-k\right|  dx  \\
 & \geq & c\int_{\Omega}|\nabla w|^2
 d^{2-k}\,dx+c\int_{\Omega}|w|^2  d^{-k}
\left|d\Delta d+1-k\right| X(d/D) dx\\
\mbox{(by (\ref{eq:sh1}))} & \geq &c\int_{\Omega}|w|^2d^{-k}X^2(d/D)dx\\
 &=&c\int_{\Omega}\frac{|u|^p}{d^p}X^2(d/D) \, dx.
\end{eqnarray*}
$\hfill //$


 \setcounter{equation}{0}
 \section{The vector field approach}


 In this section we provide an alternative proof of Improved
 Hardy Inequality, based on the appropriate use  of a suitable  vector field
 and elementary calculations. It is essential for this approach that
 all terms in the  Improved
 Hardy inequality are homogeneous with respect to $u$.
  It has the advantage that it allows us to
 compute explicit constants  for the remainder term.
  In contrast, it does not work for non-homogeneous
 inequalities. We retain the geometric assumptions
introduced in Section 2. In the theorem that follows
we consider the case $p\neq k$, while
Theorem \ref{thm:os} below concerns the degenerate case
$p=k$. The optimality of the estimates
is proved in Section 5.

Let us recall the Improved Hardy inequality, which we now  write in the form
\begin{equation}
\int_{\Omega}|\nabla u|^pdx \geq |H|^p\int_{\Omega}\frac{|u|^p}{d^p}dx
+B \int_{\Omega}\frac{|u|^p}{d^p}X^2(d/D)dx.
\label{eq:6.1}
\end{equation}
We then  have
 \begin{theorem}
Assume that condition (C) is satisfied. Then, there
exists a $D_0= D_0(k,p)>0$ such that for $D \geq D_0$,
inequality (\ref{eq:6.1}) holds true with
\[
B = \frac{p-1}{2p}|H|^{p-2}.
\]
If in addition $2 \leq p <k$, then we can take
$D_0 = \sup_{x \in \Omega} d(x,K)$.
\label{thm:ddd}
\end{theorem}

{\em Proof.}
Let $T$ be a $C^1$ vector field on $\Omega$.
For any  $u \in C_c^{\infty}(\Omega \setminus K)$ we
integrate by parts and use H\"{o}lder's inequality to obtain
 \begin{eqnarray*}
 \int_{\Omega} {\rm div} T ~  |u|^p dx & = &
 -p \int_{\Omega} ( T \cdot \nabla u ) |u|^{p-2} u dx \\
 & \leq & p \left( \int_{\Omega} | \nabla u |^p  dx \right)^{\frac{1}{p}}
 \left( \int_{\Omega} |T|^{\frac{p}{p-1}} |u|^{p} dx  \right)^{\frac{p-1}{p}} \\
 & \leq &
  \int_{\Omega}  | \nabla u |^p dx  + (p-1)
 \int_{\Omega} |T|^{\frac{p}{p-1}} |u|^{p} dx.
 \end{eqnarray*}
 We therefore arrive at
 \begin{equation}
  \int_{\Omega}  | \nabla u |^p  dx  \geq
  \int_{\Omega} ( {\rm div}  T - (p-1)  |T|^{\frac{p}{p-1}} ) |u|^{p} dx.
 \label{eq:harrell}
\end{equation}
 In view of this and (\ref{eq:6.1}), the Improved Hardy inequality will be proved once
 we establish the following pointwise inequality
 \begin{equation}
\diver  T - (p-1)  |T|^{\frac{p}{p-1}} \geq \frac{|H|^{p}}{d^p} \left(1 +\frac{p-1}{2pH^2}X^2(d/D)
\right),\quad x\in\Omega.
 \label{eq:6.2}
 \end{equation}
To proceed we now make a specific choice of $T$. We take
 \[
 T(x) = H  |H|^{p-2}   \frac{\nabla d(x)}{d^{p-1}(x)} (1 +  \frac{p-1}{pH}  X(d(x)/D)
 + a X^2(d(x)/D)).
 \]
where $a$ is a free parameter to be chosen later. In any case $a$  will be
such that the quantity $1 +  \frac{p-1}{pH} X(d/D) +aX^2(d/D)$ is positive
on $\Omega$. Note that $T(x)$
is singular at $x \in K$, but since  $u \in C_c^{\infty}(\Omega \setminus K)$
all previous calculations are legitimate. A simple computation  shows that
 \begin{eqnarray*}
  {\rm div}  T &  = & H  |H|^{p-2}  \frac{d \Delta d - (p-1) |\nabla d|^2}{d^p}
 \left(1+  \frac{p-1}{pH}  X +  aX^2(d/D)\right)  \\
 & & \,
 + H  |H|^{p-2}  \frac{ |\nabla d|^2}{d^p}\left( \frac{p-1}{pH}  X^2(d/D) +
2aX^3(d/D)\right) \\
 & \geq & H|H|^{p-2}\frac{k-p}{d^p}\left(1+\frac{p-1}{pH}X(d/D)+a X^2(d/D)\right)\\
&& +H  |H|^{p-2}  \frac{ 1}{d^p}
\left( \frac{p-1}{pH}  X^2(d/D) + 2aX^3(d/D)\right),
 \end{eqnarray*}
where in the last  inequality we used (C) and the fact that $|\nabla d|=1$.
Thus, we have
\begin{eqnarray*}
&&\hspace{-2cm}\diver T - (p-1)  |T|^{\frac{p}{p-1}} \geq \\
&\geq& H  |H|^{p-2} \frac{(k-p) (1+  \frac{p-1}{pH}  X(d/D) +  aX^2(d/D))}{d^p}+\\
&&+H |H|^{p-2}  \frac{ ( \frac{p-1}{pH}  X^2(d/D) + 2aX^3(d/D)) }{d^p}-\\
&&-(p-1)|H|^p \frac{ (1+  \frac{p-1}{pH}  X(d/D)+ aX^2(d/D))^{\frac{p}{p-1}}}{d^p}.
\end{eqnarray*}
It then follows that for (\ref{eq:6.2}) to hold,  it is enough
to establish the inequality
\begin{equation}
f(t) \geq 1  + \frac{p-1}{2pH^2}t^2, \qquad t\in [0,M],
\label{eq:6.3}
\end{equation}
where $M=M(D):=\sup_{x\in\Omega}X(d(x)/D)$ and
 \[
 f(t) := p (1+  \frac{p-1}{pH}  t +  at^2)
  + \frac{1}{H}  ( \frac{p-1}{pH}  t^2 + 2at^3)
    -(p-1)  (1+  \frac{p-1}{pH}  t +  at^2)^{\frac{p}{p-1}}.
 \]
From Taylor's formula we have that
 \begin{equation}
 f(t) = f(0) +  f'(0) t + \frac{1}{2}f''( \xi_{t}) t^2,~~~~~~~~~0\leq \xi_t \leq t \leq M.
 \label{eq:6.4}
 \end{equation}
 We have $f(0)=1$. Moreover,
 \begin{eqnarray}
 f'(t)  & = &  \frac{p-1}{H} +2ap t + \frac{2(p-1)}{pH^2}  t + \frac{6a}{H} t^2
  \nonumber \\
 & & \, -p (1+  \frac{p-1}{pH}  t +  at^2)^{\frac{1}{p-1}}
  ( \frac{p-1}{pH}  + 2at)
 ,~~~~~~~  \nonumber \\
 f''(t)  & = & 2ap + \frac{2(p-1)}{pH^2} +  \frac{12a}{H} t -
 2ap (1+  \frac{p-1}{pH}  t +  at^2)^{\frac{1}{p-1}} \label{eq:6.55}  \\
 & & -\frac{p}{p-1}  (1+  \frac{p-1}{pH}  t +  at^2)^{\frac{2-p}{p-1}}
  ( \frac{p-1}{pH}  + 2at)^2,   \nonumber \\
 f'''(t)  & = &   \frac{12a}{H} - \frac{6ap}{p-1}
  (1+  \frac{p-1}{pH}  t +  at^2)^{\frac{2-p}{p-1}} ( \frac{p-1}{pH}  + 2at)
 \nonumber  \\
 & & -\frac{p(2-p)}{(p-1)^2}  (1+  \frac{p-1}{pH}  t +  at^2)^{\frac{3-2p}{p-1}}
 ( \frac{p-1}{pH}  + 2at)^3,  \nonumber
 \end{eqnarray}
and in particular
 \begin{eqnarray}
  f'(0) & = & 0,   \nonumber  \\
  f''(0) &  =  & \frac{p-1}{pH^2},  \\
 f'''(0) &  = &  \frac{6a}{H}-\frac{(2-p)(p-1)}{p^2H^3}. \nonumber
 \label{eq:6.6}
 \end{eqnarray}
 To proceed we distinguish  various cases.

 {\bf (a)} $1<p<2 \leq k$.  In this case $H>0$.  We now choose $a$ so that
  $f'''(0)>0$, that is, $a >  \frac{2-p}{6(p-1)}>0$.
 Hence  $f''$ is an increasing function in some interval of the form $(0, M_0)$.
 Consequently, for $t \in (0,M_0)$
 \[
 f''( \xi_{t}) \geq f''(0) =  \frac{p-1}{pH^2}.
 \]
 It then follows from  (\ref{eq:6.4})
 \[
 f(t)  \geq 1 +  \frac{p-1}{2pH^2} t^2,\quad t\in [0,M_0].
 \]
Hence (\ref{eq:6.3}) has been proved in this case.

 {\bf (b)} $2 \leq p <k$. We still have $H>0$.   We now choose $a=0$.
 It is clear that $f'''(0)>0$.  Moreover, we compute
 \[
 f'''(t) =  \frac{(p-1)(p-2)}{p^2H^3}
(1+\frac{p-1}{pH}  t)^{\frac{3-2p}{p-1}} >0,\qquad{\rm for ~all}~~~t>0.
 \]
 We  then  repeat the argument of case (a),  taking $M_0=+\infty$.

 {\bf (c)} $k=1< p <2$. We now have $H<0$. We then  choose $a$ such that
 $0<a<(2-p)(p-1)/(6p^2H^2)$, so that $f'''(0)>0$ and the previous
 argument goes through.

 {\bf (d)} $p \geq 2$, $p>k$. Again $H<0$.  We now take
$a <(2-p)(p-1)/(6p^2H^2)<0$ and proceed as before.

 It is clear that we can choose an $M_0$ (small enough)
  that works simultaneously in
 all cases, and at the same time  $(1+  \frac{p-1}{pH}  X +  aX^2)>0$,
 for $0<X<M_0$. We can even estimate this $M_0$ using (\ref{eq:6.55}), if needed.
 Since $X(d/D)= -\log^{-1}(d/D)$  the condition $X \leq M_0$ is equivalent to
 $D \geq D_0 := e^{1/M_0} \sup_{x \in \Omega} d(x)$.
 The proof of the theorem is now complete.$\hfill //$

{\bf Remark} The assumption $\sup_{x \in \Omega} d(x) < +\infty$ is only needed in
order to obtain the Improved Hardy.  For the plain Hardy inequality one can 
choose the vector field  $T(x) = H  |H|^{p-2}   \frac{\nabla d(x)}{d^{p-1}(x)}$
in which case the boundedness of $d(x)$ is not required.

Clearly the usual Hardy inequality does not hold when $p=k$.
In our next result we give a substitute for Hardy inequality in that case.
The analogue of condition (C) is now
\[p=k,\qquad  d\Delta d +1-k\geq 0. \qquad\qquad {\bf\rm (C')}\]

In Theorem \ref{thm:opos} we shall prove that estimate (\ref{eq:osfp}) below is sharp.
Our result reads

\begin{theorem}
Let $p=k$ and assume that $d(\cdot)$ is bounded in $\Omega$.
If $(C')$ is satisfied then for any
$D\geq\sup_{\Omega}d(x)$ there holds
\begin{equation}
\int_{\Omega}|\nabla u|^pdx\geq \left(\frac{p-1}{p}\right)^p
\int_{\Omega}\frac{|u|^p}{d^p}X^p(d/D)dx,\quad u\in W^{1,p}_0
(\Omega\setminus K).
\label{eq:osfp}
\end{equation}
\label{thm:os}
\end{theorem}
{\em Proof. }We define the vector field
\[T(x)=\left(\frac{p-1}{p}\right)^{p-1}\frac{X^{p-1}(d(x)/D)}{d^{p-1}(x)}\nabla d(x),\quad
x\in\Omega,\]
and use inequality (\ref{eq:harrell}). We have
\begin{eqnarray*}
\diver T&=&(\frac{p-1}{p})^{p-1}d^{-p}X^{p-1}(d/D)\left((p-1)X(d/D)-p+1+d\Delta d\right)\\
&\geq&(\frac{p-1}{p})^{p-1}(p-1)d^{-p}X^p(d/D),
\end{eqnarray*}
and hence
\[\diver T -  (p-1)  |T|^{\frac{p}{p-1}}\geq (\frac{p-1}{p})^{p}d^{-p}X^p(d/D),\]
which yields (\ref{eq:osfp}).  $\hfill //$


 \setcounter{equation}{0}
 \section{Best constants for Improved Hardy}


 
In this section we will prove the optimality of the constants appearing
in the Improved Hardy Inequalities we derived in Section 4.
This will be done by deriving optimal bounds for all constants appearing
in improved Hardy inequalities of the type we consider in this work.
More precisely, recalling that $H=(k-p)/p$, we have the following:

\begin{theorem}
Let $\Omega$ be a domain in $\R^N$. $\ia$ If $2\leq k\leq N-1$ then we take $K$ to
be a piecewise smooth surface of codimension $k$ and assume $K\cap\Omega\neq 
\emptyset$; $\ib$ if $k=N$ then we take
$K=\{0\}\subset\Omega$; $\ic$ if $k=1$ then we take $K =\partial \Omega$.
Suppose that for some constants $A>0$, $B\geq 0$ and  $\gamma>0$,
the following inequality holds true for all $u\in C^{\infty}_{c}(\Omega\setminus K)$
\begin{equation}
\int_{\Omega}|\nabla u|^pdx   \geq    A \int_{\Omega}\frac{|u|^p}{d^p}dx
 +  B \int_{\Omega}\frac{|u|^p}{d^p}X^{\gamma}(d/D)dx.
\label{eq:prop}
\end{equation}
Then:

$\ia$   $A \leq |H|^p$.

$\ib$ If $A = |H|^p$,  and $B>0$, then $\gamma \geq 2$.

$\ic$ If  $A = |H|^p$  and  $\gamma = 2$, then $B  \leq 
 \frac{p-1}{2p}|H|^{p-2}$.
\label{th:best}
\end{theorem}

To prove this theorem  we  will use  a minimizing sequence for the
Improved Hardy inequality. Without any loss of generality we may
assume that $0\in K\cap\Omega$.
All our analysis will be local, say, in a fixed ball of radius  $\delta$
(denoted by $B_{\delta}$) centered 
at the origin, for some fixed small $\delta$.
We next introduce the function
\begin{equation}
w_{\epsilon}(x) = d^{-H+\epsilon}(x) X^{-\theta}(d(x)/D),~~~~~~~1/p <\theta <2/p,
\label{eq:test1}
\end{equation}
where $~k=1,2,\ldots,N$ and
 $D=\sup_{x \in \Omega} \dist(x,K)$, as usual.  In order to localize it 
we also define a suitable non-negative  test function $\phi \in C_c^{2}(B_{\delta})$,
such that $\phi(x)=1$ for  $x \in B_{\delta/2}$.
 We then set
\begin{equation}
U_{\epsilon}(x) = \phi(x) w_{\epsilon}(x),~~~
~~~~\supp U_{\epsilon} \subset  B_{\delta}.
\label{eq:test2}
\end{equation}

The proof we present works for any $k=1,2,\ldots,N$. We note however that
for $k=N$ (distance from a point) the subsequent calculations are substantially
simplified, whereas for $k=1$ (distance from the boundary) one should
replace $B_{\delta}$ by
$B_{\delta} \cap \Omega$. This last  change entails some minor modifications,
the arguments otherwise being  the same.

Throughout the rest of this Section we denote by $C$, $c(p)$ etc
various positive constants, not necessarily the same in each occurrence,
which may depend on $\delta$, $p$ or $k$ but are independent of $\epsilon$.

We begin by presenting some lemmas that contain all technical estimates that
we need for the proof of the theorem. For $\beta \in \R$ and $\epsilon>0$ small we define 
\begin{equation}
J_{\beta}(\epsilon)=\int_{\Omega}\phi^pd^{-k+\epsilon p} X^{-\beta}(d/D)dx.
\label{eq:j}
\end{equation}

\begin{lemma}
For $\epsilon$ small there holds                   
\begin{eqnarray*}
\ia && c\epsilon^{-1-\beta}\leq J_{\beta}(\epsilon)\leq
C\epsilon^{-1-\beta},\quad\quad\mbox{$\beta>-1$;}\\
\ib && J_{\beta}(\epsilon)=\frac{p\epsilon}{\beta+1}J_{\beta+1}(\epsilon) +O_{\epsilon}(1),\quad
\mbox{$\beta>-1$;}\\
\ic && J_{\beta}(\epsilon) = O_{\epsilon}(1),\quad\quad\mbox{$\beta<-1$.}
\end{eqnarray*}
\label{lem:j}
\end{lemma}
{\em Proof. } Since $|\nabla d|=1$ we have
\[ J_{\beta}(\epsilon)=\int_0^{\delta}\int_{d=r}\phi^pr^{-k+\epsilon p}X^{-\beta}(r/D)dS\,dr.\]
Hence using the fact $0\leq\phi\leq 1$ and $\int_{\{d=r\}\cap B_{\delta}} dS <cr^{k-1}$, we obtain
\begin{eqnarray*}
J_{\beta}(\epsilon)  \leq   c \int_0^{\delta} r^{-1+\epsilon p}X^{-\beta}(r/D)dr.  
\end{eqnarray*}
Recalling (\ref{eq:int})
we see that for $\beta <-1$ the integral above has a finite limit as $\epsilon \rightarrow 0$,
hence (iii) follows. To show (i) we  use
 the change of variables $r=Ds^{1/\epsilon}$  to obtain that
\[
J_{\beta}(\epsilon)  \leq  \epsilon^{-1-\beta} D^{\epsilon p} \int_0^{(\delta/D)^{\epsilon}}
s^{p-1} X^{-\beta}(s) ds,
\]
and the upper estimate of (i) follows.  For the lower
estimate we use the fact that $\phi=1$ for $d\leq\delta /2$ and argue similarly.

To prove (ii) we recall (\ref{eq:de}) to write
\[(\beta+1)J_{\beta}(\epsilon)=-\int_{\Omega}\phi^pd^{1-k+\epsilon p}\nabla d
\cdot\nabla X^{-\beta-1}(d/D)dx.\]
We now perform an integration by parts and note that no boundary
terms appear. Indeed, if $k=1$ then the factor $d^{1-k+\epsilon p}=d^{\epsilon p}$
guarantees that the integrand vanishes on $K$. If $2\leq k\leq N$
then we approximate $\Omega$ by $\Omega_{\eta}:=\{x\in\Omega : d(x)>\eta\}$,
$\eta>0$ small. This yields the boundary term
\[-\int_{d=\eta}\phi^p d^{1-k+\epsilon p}X^{-\beta-1}(d/D)\nabla d\cdot
\vec{n}\, dS\]
which vanishes as $\eta\to 0$. Hence in any case we have
\begin{eqnarray*}
(\beta+1)J_{\beta}(\epsilon) & = & \int_{\Omega}\diver
(\phi^p d^{1-k+\epsilon p}\nabla d) X^{-\beta-1}(d/D)dx\\
&=&p\int_{\Omega}\phi^{p-1}d^{1-k+\epsilon p}X^{-\beta-1}(d/D)
\nabla\phi\cdot\nabla d\; dx +\\
&&+(1-k+\epsilon p)\int_{\Omega}d^{-k+\epsilon p}X^{-\beta-1}(d/D)dx+\\
&&+\int_{\Omega}\phi^pd^{1-k+\epsilon} \Delta d X^{-\beta-1}(d/D)dx.
\end{eqnarray*}
The first integral is of order $O_{\epsilon}(1)$ by an application
of (i). The other two integrals combine to give
\begin{equation}
\epsilon p J_{\beta+1}(\epsilon) +\int_{\Omega}\phi^p
d^{-k+\epsilon p}X^{-\beta-1}(d/D)(d\Delta d+1-k)dx.
\label{eq:nel}
\end{equation}
But it is a direct consequence of \cite[Theorem 3.2]{AS} that
\begin{equation}
d \Delta d + 1-k = O(d)~~~~~~{\rm  as} ~d(x) \rightarrow 0;
\label{mete}
\end{equation}
this implies that the integral in (\ref{eq:nel}) is of order $O_{\epsilon}(1)$,
and the result follows. $\hfill //$

We next estimate the quantity 
\[
I[ U_{\epsilon}] = \int_{\Omega}|\nabla U_{\epsilon} |^pdx  -
|H|^p  \int_{\Omega}\frac{| U_{\epsilon}|^p}{d^p}dx.
\]
\begin{lemma}
As $\epsilon \rightarrow 0$, there holds
\begin{eqnarray}
\ia &&
I[U_{\epsilon}]\leq \frac{\theta (p-1)}{2}
|H|^{p-2} J_{p\theta -2}(\epsilon) +O_{\epsilon }(1); 
\label{eq:l5.2}\\
\ib &&\int_{B_{\delta}}|\nabla U_{\epsilon}|^pdx  \leq 
 |H|^p J_{p\theta}(\epsilon)+O_{\epsilon }(\epsilon^{1-p\theta}).
\label{eq:l5.1}
\end{eqnarray}
\label{lem:5.2}
\end{lemma}
{\em Proof.} We have $\nabla U_{\epsilon}= 
 \phi\,\nabla w_{\epsilon}+\nabla\phi\, w_{\epsilon} $ and
hence, using the elementary inequality
\begin{equation}
|a+b|^p  \leq |a|^p + c_p (|a|^{p-1} |b| + |b|^p), ~~~a,b \in \R^N, ~~~~~p>1,
\label{eq:tv}
\end{equation}
we obtain
\begin{eqnarray}
\int_{\Omega}|\nabla U_{\epsilon}|^pdx 
 &  \leq  & \int_{B_{\delta}} \phi^{p} 
 d^{-k +\epsilon p} X^{-p \theta}(d/D)
\left| H- (\epsilon  -\theta X(d/D)) \right|^p dx \nonumber   \\
& & \!\! +c_p \,  \int_{B_{\delta}}  | \nabla \phi| |\phi|^{p-1}
 |\nabla  w_{\epsilon} |^{p-1}
|w_{\epsilon}|\, dx  +c_p \int_{B_{\delta}} | \nabla \phi|^p | w_{\epsilon}|^p\, dx
\nonumber  \\
& := & I_A + I_2 + I_3.  \label{eq:l5.1.2}
\end{eqnarray}
We claim that
\begin{equation}
I_2,~I_3  = O_{\epsilon}(1), ~~~~~~~~~~{\rm as} ~~\epsilon \rightarrow 0.
\label{eq:i23}
\end{equation}
Let us give the proof for $I_2$. Using the definition of $w_{\epsilon}$ and
the fact that $\phi$ is a nice function we get
\[
I_2 \leq c\int_{B_{\delta}} d^{1-k +\epsilon p} X^{-p \theta}(d/D)
\left| H- (\epsilon  -\theta X(d/D)) \right|^{p-1}  dx.
\] 
Since  $(\epsilon  -\theta X(d/D))$ is small compared to $H$ we have
\[
I_2 \leq c \int_{B_{\delta}} d^{1-k +\epsilon p} X^{-p \theta}(d/D) dx.
\]
The integral in the right hand side has a finite limit as $\epsilon\to 0$ by (i)
of Lemma \ref{lem:j}. The integral $I_3$ is treated in the same way.

From (\ref{eq:l5.1.2}), (\ref{eq:i23}) and the definition of $J_{\beta}$ we
easily obtain 
\begin{eqnarray}
I[ U_{\epsilon}] &= &\int_{B_{\delta}}|\nabla U_{\epsilon} |^pdx - |H|^{p} J_{p \theta} \nonumber \\
& \leq & I_A - |H|^{p} J_{p \theta} + O_{\epsilon}(1)\label{eq:l5.21}  \\
& = &   I_1 + O_{\epsilon }(1),  \nonumber 
\end{eqnarray}
where
\[
I_1 :=
\int_{B_{\delta}}\phi^p 
 d^{-k +\epsilon p} X^{-p \theta}(d/D)
 \Bigl(\left| H- \bigl(\epsilon  -\theta  X(d/D)\bigr) \right|^p
 -|H|^p\Bigr)dx. 
\]
We proceed by estimating $I_1$. Since $ \eta :=(\epsilon  -\theta X(d/D))$
is  small compared to $H$ we may use Taylor's expansion  to obtain
\[
|H - \eta|^p - |H|^p \leq - p |H|^{p-2} H \eta + \frac{1}{2}p(p-1) |H|^{p-2} \eta^2
+ C |\eta|^3.
\]
Using this inequality we can bound $I_1$ by
\begin{equation}
I_1 \leq I_{11} + I_{12} + I_{13},
\label{eq:l5.22}
\end{equation}
where
\begin{eqnarray*}
I_{11}&  = &  - p |H|^{p-2} H 
 \int_{B_{\delta}} \phi^{p} 
 d^{-k +\epsilon p} X^{-p \theta}(d/D)
  (\epsilon  -\theta X(d/D))dx, \\
I_{12}&  = &  \frac{1}{2}p(p-1) |H|^{p-2}
 \int_{B_{\delta}} \phi^{p} 
 d^{-k +\epsilon p} X^{-p \theta}(d/D)
  \bigl( \epsilon -\theta X(d/D)\bigr)^2dx, \\
I_{13}&  = & C
 \int_{B_{\delta}} \phi^{p} 
 d^{-k +\epsilon p} X^{-p \theta}(d/D)
  |\epsilon  -\theta X(d/D)|^3 dx.
\end{eqnarray*}
We shall prove that
\begin{equation}
I_{11},I_{13}=O_{\epsilon}(1),\quad \epsilon\to 0.
\label{eq:l5.23}
\end{equation}
Indeed, an application of part (ii) of Lemma \ref{lem:j} (with $\beta=-1+p\theta$)
gives $I_{11}=O_{\epsilon}(1)$ for small $\epsilon>0$.
Concerning $I_{13}$ we clearly  have
\[I_{13} \leq  c \epsilon^3 J_{p\theta} +cJ_{p\theta -3}.\]
It follows from Lemma \ref{lem:j} parts (i) and (iii), 
and the fact that $1<p\theta<2$
 that  $I_{13}=O_{\epsilon}(1)$.

To calculate the term $I_{12}$ we first expand the square in the integrand
and then apply (ii) of Lemma \ref{lem:j} twice (with $\beta = -1+p\theta$ the
first time and $\beta = -2+p\theta>-1$ the second time) to conclude that
\begin{equation}
I_{12} = \frac{\theta (p-1)}{2} |H|^{p-2} 
\int_{B_{\delta}} \phi^p
 d^{-k +\epsilon p} (-\log d/D)^{p \theta-2} dx +  O_{\epsilon }(1),
\qquad\epsilon \rightarrow 0.
\label{eq:l5.24}
\end{equation}
From (\ref{eq:l5.21}), (\ref{eq:l5.22}), (\ref{eq:l5.23}) and 
(\ref{eq:l5.24}) we conclude 
(\ref{eq:l5.2}). The second inequality of the lemma
follows from the first  equality in (\ref{eq:l5.21}), estimate (\ref{eq:l5.2}) and Lemma \ref{lem:j}
 $\hfill //$

We are now ready to give the proof of Theorem \ref{th:best}

{\em Proof of Theorem \ref{th:best}.} It follows directly from
part (i) of Lemma \ref{lem:j} that for any $\gamma\in\R$ there holds

\begin{equation}
R_{\gamma}[U_{\epsilon}]:= \int_{\Omega}\frac{|U_{\epsilon}|^p}{d^p}X^{\gamma}(d/D)dx =J_{p\theta-\gamma}
(\epsilon).
\label{eq:pr2}
\end{equation}

{\bf (i)} Since inequality (\ref{eq:prop}) holds for every 
 $u\in W^{p}_{0}(\Omega\setminus K)$ we have

\begin{eqnarray*}
A &\leq& \frac{\int_{B_{\delta}}|\nabla U_{\epsilon}|^pdx }{R_0[U_{\epsilon} ]}\\
\mbox{(by (\ref{eq:l5.1})) }&\leq&\frac{|H|^p(1+O_{\epsilon}(\epsilon))
J_{p\theta}(\epsilon)+O_{\epsilon}(1)}{J_{p\theta}(\epsilon)};
\end{eqnarray*}
letting $\epsilon\to 0$ and recalling that $J_{p\theta}(\epsilon)
\to\infty$ we conclude that $A\leq |H|^p$.

{\bf (ii)} Let  $A=|H|^p$. Assuming that $\gamma<2$ we will reach
a contradiction. Since $p\theta-\gamma>-1$
arguing as in (i) we have that
\begin{eqnarray*}
0< B  &\leq& \frac{I[U_{\epsilon}]}{R_{\gamma}[U_{\epsilon}]}\\
\mbox{(by (\ref{eq:l5.2}) and Lemma 5.2 (i) )}&\leq&
\frac{c\epsilon^{-p\theta +1}}{c\epsilon^{-1-p\theta+\gamma}}\\
&=&c\epsilon^{2-\gamma}\to 0,\qquad\mbox{as $\epsilon\downarrow 0$,}
\end{eqnarray*}
which is a contradiction. Hence $\gamma \geq 2$.

{\bf (iii)} If $A=|H|^p$ and  $\gamma=2$ then
\begin{eqnarray*}
B&\leq&\frac{I[U_{\epsilon}]}{R_2[U_{\epsilon}]}\\
\mbox{(by (\ref{eq:l5.2}))}&\leq&\frac{\frac{1}{2}\theta (p-1) |H|^{p-2}
J_{p\theta -2}(\epsilon) +O_{\epsilon }(1)}{J_{p\theta-2}}.
\end{eqnarray*}
The assumption $\theta>1/p$ implies $J_{p\theta}\to\infty$ as $\epsilon\to 0$
by (i) of Lemma \ref{lem:j}. Hence $B\leq \frac{\theta (p-1)}{2} |H|^{p-2}$;
letting $\theta\to 1/p$ concludes the proof.$\hfill //$

We close this section proving the optimality of the estimate
in Theorem \ref{thm:os}.

\begin{theorem}
Let $\Omega$ be a domain in $\R^N$. $\ia$ If $2\leq k\leq N-1$ then we take $K$ to
be a piecewise smooth surface of codimension $k$ and assume $K\cap\Omega\neq 
\emptyset$; $\ib$ if $k=N$ then we take
$K=\{0\} \subset \Omega$.
Suppose that $p=k$ and that for some constants $B\geq 0$ and  $\gamma>0$
the following inequality holds true for all $u\in C^{\infty}_{c}(\Omega\setminus K)$
\begin{equation}
\int_{\Omega}|\nabla u|^pdx   \geq    B \int_{\Omega}\frac{|u|^p}{d^p}X^{\gamma}(d/D)dx.
\label{eq:prop1}
\end{equation}
We then have:

$\ia$ If $B>0$, then $\gamma \geq p$.

$\ib$ If $\gamma = p$, then $B  \leq  \bigl(\frac{p-1}{p}\bigr)^p$.
\label{thm:opos}
\end{theorem}

{\em Proof.} The proof uses an argument similar to that of Theorem \ref{th:best}.
Without any loss of generality we assume that $0\in K\cap \Omega$ if $2\leq k\leq N$
and $0\in\partial\Omega =K$ if $k=1$.
As in the last theorem we let $\phi$ be a non-negative,
smooth cut-off function supported in $B_{\delta}=\{|x|<\delta\}$ and equal
to one on $B_{\delta/2}$. For any $\epsilon>0$ and $\theta>(p-1)/p$ define
$w_{\epsilon}=d^{\epsilon}(-\log d)^{\theta}$ and
$U_{\epsilon}=\phi d^{\epsilon}(-\log d)^{\theta}$.
Using (\ref{eq:tv}) we have
\begin{eqnarray*}
\int_{\Omega}|\nabla U_{\epsilon}|^pdx&\leq &\int_{B_{\delta}} \phi^{p} 
|\nabla w_{\epsilon}|^pdx+
c_p\int_{B_{\delta}}\phi^{p-1}|\nabla w_{\epsilon}|^{p-1}w_{\epsilon}|\nabla\phi|dx+
c_p\int_{B_{\delta}}w_{\epsilon}^p|\nabla\phi|^p\\
&=:&I_A+I_2+I_3.
\end{eqnarray*}
Arguing as in the proof of the previous theorem we see that
\[ I_2,I_3= O_{\epsilon}(1),\qquad \epsilon\to 0.\]
Denoting by $c_i^p$ the coefficients of the binomial expansion we have 
\begin{eqnarray*}
|\nabla w_{\epsilon}|^p&=&d^{-k+\epsilon p}X^{-p\theta}(d/D)|\epsilon-\theta X(d/D)|^p\\
&\leq&d^{-k+\epsilon p}\theta X^{-p\theta}(d/D)(\epsilon +\theta X(d/D))^p\\
&=&d^{-k+\epsilon p}X^{-p\theta}(d/D)\sum_{i=0}^p c^p_i\epsilon^{p-i}\theta^{i}X^{i}(d/D),
\end{eqnarray*}
and hence
\[I_A\leq \sum_{i=0}^pc^p_i\epsilon^{p-i}\theta^{i}J_{p\theta -i}(\epsilon),\]
where the functions
$J_{\beta}(\epsilon)=\int_{\Omega}\phi^pd^{-p+\epsilon p}X^{-\beta}(d/D)$
are as in (\ref{eq:j}). Now it follows from (ii) of Lemma \ref{lem:j} and a simple induction
argument that
\[\epsilon^{p-i}J_{p\theta -i}=(\theta-\frac{i}{p})(\theta-\frac{i+1}{p})
\ldots(\theta-\frac{p-1}{p})J_{p\theta -p}+O_{\epsilon}(1),\quad i=0,\ldots,p-1.\]
But the fact that $\theta>(p-1)/p$ implies that $J_{p\theta -p}(\epsilon)\to +\infty$ as
$\epsilon\to 0$, by (i) of Lemma \ref{lem:j}. It follows that
\begin{eqnarray*}
B&\leq&\limsup_{\epsilon\to 0}\frac{\int_{\Omega}|\nabla U_{\epsilon}|^pdx}{\int_{\Omega}
\frac{U_{\epsilon}^p}{d^p}X^p(d/D)dx}\\
&\leq&\limsup_{\epsilon\to 0}\frac{\left(\theta^p
\sum_{i=0}^{p-1}c^p_i\theta^{i}(\theta-\frac{i}{p})\ldots(\theta-\frac{p-1}{p})\right)
J_{p\theta -p}(\epsilon)+O_{\epsilon}(1)}
 {J_{p\theta -p}(\epsilon)}\\
&=&\theta^p+\sum_{i=0}^{p-1}c^p_i\theta^{i}(\theta-\frac{i}{p})\ldots(\theta-\frac{p-1}{p})
\end{eqnarray*}
This last expression converges to $(\frac{p-1}{p})^p$ as $\theta\to (p-1)/p$ ; this completes
the proof.$\hfill //$




 \setcounter{equation}{0}
 \section{Non-homogeneous Improved Hardy inequalities}
 As an application of the results in Section 3
 we first prove Theorem B, the improved Hardy-Poincar\'{e} inequality.

 {\em Proof of Theorem B.} We shall prove that
\begin{equation}
 I[u]\geq
 c \left(\int_{\Omega}|\nabla u|^qd^{k(-1+q/p)}X^{\beta}(d/D)dx\right)^{p/q},
 \quad u\in W^{1,p}_0(\Omega \setminus K).
 \label{eq:hp2}
 \end{equation}
Letting $v=ud^{(k-p)/p}$ we have
 \begin{equation}
 \int_{\Omega}|\nabla u|^qd^{k(-1+q/p)}X^{\beta}dx
 \leq c(q,k)\left(\int_{\Omega}|\nabla v|^q
 d^{q-k}X^{\beta}dx +\int_{\Omega}|v|^qd^{-k}X^{\beta}dx\right).
 \label{eq:com}
 \end{equation}
 To proceed we will estimate the two integrals in the right hand side of
 (\ref{eq:com}).

 We first consider the case $1<p<2$. Applying
 H\"{o}lder's inequality we obtain
 \[
 \int_{\Omega}|\nabla v|^q
 d^{q-k}X^{\beta}dx \leq
 \left(\int_{\Omega}|\nabla v|^pd^{p-k}X^{2-p}dx \right)^{q/p}
 \left(\int_{\Omega}  d^{-k} X^{\theta} dx \right)^{(p-q)/p},
 \]
 with $\theta = (\beta p-2q+qp)/(p-q)$. We next show that the last integral above
 is finite.  The integrand has a singularity as $d \rightarrow 0$. However for
 $d$ near zero the integral behaves like
 \[
 \int_{0}^{\epsilon} \int_{d=t} \frac{d^{-k} X^{\theta}}{|\nabla d|} dS dt
 \leq c
 \int_{0}^{\epsilon} \int_{d=t} d^{-k} X^{\theta} dS dt
 \leq c  \int_{0}^{\epsilon} t^{-1} X^{\theta}(t) dt.
 \]
 The last integral is finite iff $\theta>1$ (cf (\ref{eq:int})),
 a condition that is easily seen to
 be satisfied under our assumptions on $p$, $q$, $\beta$. Hence we end up with
 \begin{equation}
 \left(\int_{\Omega}|\nabla v|^q
 d^{q-k}X^{\beta}dx \right)^{p/q}  \leq c
 \int_{\Omega}|\nabla v|^pd^{p-k}X^{2-p}dx
 \leq c I[u],
 \label{eq:ho1}
 \end{equation}
where in the last inequality we used (\ref{eq:1p2}). Applying in a similar
fashion H\"{o}lder's inequality and then the Improved Hardy inequality (\ref{eq:1p2}),
we  estimate the last integral in (\ref{eq:com})
 \begin{equation}
 \left(\int_{\Omega}|v|^qd^{-k}X^{\beta}dx  \right)^{p/q}
 \leq c \int_{\Omega}|v|^pd^{-k}X^2dx =
 \int_{\Omega} \frac{|u|^p}{d^{k}}X^2  dx
 \leq c I[u]
 \label{eq:ho2}
 \end{equation}
 and (\ref{eq:hp2}) follows.

 Consider now the case  $p \geq 2$. The proof is quite similar. In particular
 estimate (\ref{eq:ho2}) remains valid, whereas the analogue of (\ref{eq:ho1})
 is
 \[
 \left(\int_{\Omega}|\nabla v|^q
 d^{q-k}X^{\beta}dx \right)^{p/q}
 \leq c
 \int_{\Omega}|\nabla v|^p
 d^{p-k} dx  \leq c I[u],
 \]
where in the last inequality we used (\ref{eq:aa11}) and condition (C).
The proof of (\ref{eq:hp2}) is now  complete.

To prove the sharpness of the estimate we consider the functions
$U_{\epsilon}$ of Section 5 (see (\ref{eq:test2})). We have already seen that they satisfy
$I[U_{\epsilon}]\leq c\epsilon^{1-p\theta}$. Moreover, simple
calculations show that for $\beta>0$ we have
\[\int_{\Omega}|\nabla U_{\epsilon}|^qd^{k(-1+q/p)}X^{\beta}(d/D)dx
\geq c\epsilon^{\beta-\theta q-1},\]
for all $\epsilon>0$ small enough. Hence (\ref{eq:hp2}) cannot
be true if $\beta< 1+q/p$.$\hfill //$


 We now turn our attention to improved Hardy-Sobolev
 inequalities. By this we mean lower estimates on $I[u]$
 in terms of weighted $L^q$ norms of the function $u$, $q>p$.
 It will be seen that a there is a difference in
 the form the estimates take, depending on
 whether $k=N$ or $k< N$. We first consider the case
 of affine $K$, $K\equiv \R^{N-k}$, and take $\Omega=\R^N$.
 More precisely, we write points in $\R^N$
 as $x=(y,z)$, $y\in\R^{N-k}$, $z\in\R^k$. Under
 this representation we take
 \[ K=\{(y,0)\,|\, y\in\R^{N-k}\}\]
so that $d(y,z)=|z|$.

Our next two propositions yield Theorem C.
 \begin{proposition}
Assume that $k< N$ and that condition (C) is satisfied.
Then for any $2\leq p<N$ and any
 $p<q\leq Np/(N-p)$ there exists $c>0$ such that
 \begin{equation}
 I[u]\geq c\left(\int_{\Omega}|u|^q|z|^{-q-N+Nq/p}dx
 \right)^{p/q},\quad u\in W^{1,p}_0(\R^N\setminus K).
 \label{eq:m1}
 \end{equation}
 \label{thm:mm1}
 \end{proposition}
 {\em Proof. }Let $v(y,z)=u(y,z)|z|^{(k-p)/p}$. It follows from
 (\ref{eq:aa11}) and condition (C) that
 \[I[u]\geq c\int_{\R^N}|\nabla v|^p |z|^{p-k}dzdy.\]
 Moreover, Corollary 2 Section 2.1.6 of \cite{M} gives
 \begin{equation}
 \int_{\R^N}|\nabla v|^p |z|^{p-k}dzdy\geq
 c\left(\int_{\R^N}|v|^q|z|^{-N+(N-k)q/p}dzdy
 \right)^{p/q}\hspace{-10pt},v\in C^{\infty}_c(\R^N\setminus K).
 \label{eq:maz}
 \end{equation}
 Combining the last two inequalities we obtain (\ref{eq:m1}).$\hfill //$

Estimate (\ref{eq:maz}) is not valid when $k=N$ and $K$ reduces
to the single point $0\in \Omega$. Indeed, it is remarkable
that (\ref{eq:m1}) fails in this case. In our next proposition we use
decreasing rearrangement techniques to obtain a modified version
of Proposition \ref{thm:mm1} which involves a logarithmic correction
$X^{1+q/p}$ in the right hand side; we then show that the exponent $1+q/p$ is
optimal.

 \begin{proposition}
 Let  $1<p<N$.
 Let $\Omega\subset\R^N$ be bounded
 containing the origin and $D> \sup_{x \in \Omega}d(x)$.
 For any $p<q<Np/(N-p)$.
 there exists $c=c(p,q,N,\Omega)>0$ such that
 \begin{equation}
 I[u]\geq  c\left(\int_{\Omega}|u|^q |x|^{-q-N+Nq/p}X^{1+q/p}(|x|/D)dx
 \right)^{p/q},~~~
 u\in W^{1,p}_0(\Omega).
 \label{eq:1}
 \end{equation}
 Moreover one cannot replace $X^{1+q/p}$ by a lower power
 of $X$.
 \label{thm:orsob}
 \end{proposition}
 {\em Proof.} We may assume that $D=1$.
 Let $u\in C^{\infty}_c(\Omega)$ be given and let
 $u^*$ denote its radially symmetric decreasing
 rearrangement on the ball $\Omega^*$ having the same
 volume as $\Omega$ and centered at the origin.
 It is a standard property of decreasing rearrangements that
 \[I[u]\geq I[u^*].\]
 Define
 \[f(r)=r^{-q-N+Nq/p}X^{1+q/p}(r);\]
 and note that this decreases near $r=0$.
 Let $f^*:\Omega^*\to [0,\infty)$
 be the symmetric decreasing rearrangement
 of $f(|\cdot|):\Omega\to [0,\infty)$.
 Using Lemma 4.1 of \cite{FT} one sees that
 $f^*(r)\leq f(r)$, near $r=0$.
 Hence, using also the standard
 relations $\int_{\Omega}fg\leq \int_{\Omega^*}f^*g^*$ ($f,g\geq 0$),
 and $(|u|^q)^*=|u^*|^q$, we conclude that
 it is enough to establish (\ref{eq:1}) in the case
 where $\Omega$ is the unit ball and $u=u(r)$
 is a radially symmetric decreasing function of $r=|x|$.

 Assume first that  $1<p<2$  and set $v(r)=u(r)r^{(N-p)/p}$.
 Using first (\ref{eq:1p2}) (with $d=r$, $k=N$) and then
 Lemma \ref{lem:4500} (with $\alpha = 2-p$) we have
 \begin{eqnarray*}
 I[u]&\geq& c\int_0^1 |v'|^pr^{p-1}X^{2-p}dr\\
 &\geq&c\left(\int_0^1 |v|^qr^{-1}X^{1+q/p}dr\right)^{p/q}\\
 &=&c\left(\int_0^1 |u|^qr^{-q-1+Nq/p}X^{1+q/p}dr\right)^{p/q}\\
 &=&c\left(\int_{\Omega} |u|^q|x|^{-q-N+Nq/p}X^{1+q/p}(|x|)dx\right)^{p/q}.
 \end{eqnarray*}

 Suppose now that $p \geq 2$. Let $w=|v|^{p/2}$ with $v$ as above.
 Using first (\ref{eq:a71}) -- with $d=r$, $k=N$ -- and then
 Lemma \ref{lem:4500} (see Appendix) -- with $\alpha = 0$ and $2q/p$ in the place
 of $q$ -- we have
 \begin{eqnarray*}
 I[u]&\geq& c\int_0^1 |v'|^2 |v|^{p-2} rdr\\
 & = &c\int_0^1|w'|^2rdr\\
 &\geq&c\left(\int_0^1|w|^{2q/p}r^{-1}X^{1+q/p}dr\right)^{p/q}\\
 &=&c\left(\int_{\Omega}|u|^q|x|^{-q-N+Nq/p}X^{1+q/p}(|x|)dx\right)^{p/q}.
 \end{eqnarray*}
 as required.

To prove that the exponent $1+q/p$ is optimal we consider once
again the functions
$U_{\epsilon}$ of Section 5, $U_{\epsilon}(x)=\phi(x)
|x|^{\epsilon-(N-p)/p}X^{-\theta}(|x|/D)$, $\epsilon>0$, $\theta>1/p$,
$\phi$ a cut-off. An argument similar to that used in Section 5
shows the optimality of the exponent $1+q/p$.
We omit the details. $\hfill //$

Given that the estimate $I[u]\geq c\|u\|_{L^{Np/(N-p)}}^p$ is not valid,
one may ask whether the next best thing is true, i.e. whether
\[I[u]\geq c\|u\|_{L^{Np/(N-p),\infty}},\quad u\in W^{1,p}(\Omega),\]
where in the right hand side we have the weak $L^{Np/(N-p)}$ norm of $u$,
\[
  \|u\|_{L^{q, \infty}} =
 \sup_{E \subset \Omega} |E|^{\frac{1}{q}-1} \int_{E} |u| dx,\qquad 1<q<\infty.
 \]
This question was risen in a different context in \cite{BL} where Improved Sobolev Inequalities
are considered. In that paper the authors obtain lower estimates on
\[J[u]:=\int_{\Omega}|\nabla u|^2dx-c_*\left(\int_{\Omega}|u|^{2N/(N-2)}dx\right)^{(N-2)/N},
\quad u\in W^{1,2}_0(\Omega),\]
$c_*$ being the best Sobolev constant.
It is shown in \cite{BL} that $J[u]\geq c\|u\|_q^2$, $c>0$, when $q<N/(N-2)$. This of course
fails
for the critical value $q=N/(N-2)$, but it is shown instead that
\[J[u]\geq c\|u\|_{L^{N/(N-2),\infty}},\quad u\in W^{1,2}_0(\Omega).\]
In our case there is no room even for such a weak norm as the following proposition
shows.

 \begin{proposition}
There exists no $c>0$ such that
\begin{equation}
I[u]\geq c\|u\|_{L^{\frac{Np}{N-p},\infty}},\quad u\in W^{1,p}_0(\Omega).
\label{eq:ort}
\end{equation}
\label{prop:weak}
 \end{proposition}
{\em Proof.} Let $U_{\epsilon}$ be the functions introduced in
Section 5 and assume that $1/(p-1)<\theta<1/p$. We claim  that
\begin{equation}
  \|U_{\epsilon}\|_{{\frac{Np}{N-p}, \infty}(B_1)}  \geq c \epsilon^{-\theta},
 \quad ~~~~~~
 \mbox{ small $\epsilon>0$.}
 \label{eq:wlp}
\end{equation}
 Let $B_{\rho}$ denote the ball of radius $\rho$ centered at the origin. We then  have that
 \[
 \|U_{\epsilon}\|_{L^{ \frac{Np}{N-p}, \infty}(B_1)}
 \geq
 \sup_{0<\rho<1} |B_{\rho}|^{- \frac{Np-N+p}{Np}}  \int_{B_{\rho}} |u_{\epsilon}| dx
 \]
 It is immediate that $|B_{\rho}| = C \rho^{ \frac{Np-N+p}{p}}$. On the other hand
 using the explicit value of $u_{\epsilon}$  and integrating  once by parts we get
 \begin{eqnarray*}
  \int_{B_{\rho}} |U_{\epsilon}| dx  & = &
 C \int_{0}^{\rho} r^{N-\frac{N}{p} + \epsilon} (-\log r )^{\theta} dr \\
 & = & C( \rho^{N-\frac{N}{p} + \epsilon+1} (-\log r )^{\theta}
 + \int_{0}^{\rho} r^{N-\frac{N}{p} + \epsilon} (-\log r )^{\theta-1} dr )\\
 & \geq & C  \rho^{N-\frac{N}{p} + \epsilon+1} (-\log r )^{\theta}.
 \end{eqnarray*}
 Hence, we arrive at
 \[
 \|U_{\epsilon}\|_{L^{\frac{Np}{N-p}, \infty}(B_1)}
 \geq C
 \sup_{0<\rho<1} \rho^{\epsilon} (-\log \rho)^{\theta}.
 \]
 It is easy to check that $\max_{0<\rho<1}  \rho^{\epsilon} (-\log \rho)^{\theta}=
 (\theta/e)^{\theta} \epsilon^{-\theta}$ and (\ref{eq:wlp}) follows.

 On the other hand, we have seen in Section 5 that for small $\epsilon$
 \begin{equation}
 I[U_{\epsilon}]  \leq c \epsilon^{1-p\theta}.
 \label{eq:wke}
 \end{equation}
Combining (\ref{eq:wlp}) and (\ref{eq:wke}) we obtain the stated result.$\hfill //$

We close this section presenting an improved
Hardy-Sobolev inequality that is valid for all
$k\leq N$ without assuming that $K$ is affine.
The estimate obtained is weaker than that of Theorem C.

 \begin{theorem}
Let $k\leq N$ and $1<p<N$.
 Let  $D > \sup_{x \in\Omega} d(x,K)$, and assume condition (C)
is satisfied. For any $p<q\leq Np/(N-p)$ there exists
 $c=c(p,q,D,\Omega,K)>0$ such that
 \begin{equation}
 I[u] \,  \geq c \,
 \left(\int_{\Omega}|u|^qd^{-q-N+Nq/p}X^{2q/p}(d/D)dx
 \right)^{p/q},
 \label{eq:des}
 \end{equation}
 for all $u\in W^{1,p}_0(\Omega \setminus K)$.
 \label{thm:sobboun}
 \end{theorem}
 {\em Proof:}  We may assume as usual that $D= 1$. Once more we set
 $u = v d^{-H}$, $H=(k-p)/p$. From Lemma \ref{lem:sobboun} (see Appendix) 
we have
 \begin{eqnarray}
 &&\hspace{-1.5cm}\int_{\Omega}|u|^qd^{-q-N+Nq/p}X^{2q/p}dx=  \nonumber \\
 &=&\int_{\Omega}|v|^qd^{-N+(N-k)q/p}X^{2q/p}dx   \nonumber \\
 &\leq&c\left(\int_{\Omega}|\nabla v|^pd^{p-k}X^2dx +
 \int_{\Omega}|v|^pd^{-k}X^2dx\right)^{q/p}.
 \label{s1}
 \end{eqnarray}
 The last integral in (\ref{eq:s1}) is easily  estimated  by the Improved Hardy
 inequality
 \begin{equation}
 \int_{\Omega}|v|^pd^{-k}X^2dx =
  \int_{\Omega}
 \frac{|u|^p}{d^p}X^2 \, dx
 \leq  c I[u].
 \label{s2}
 \end{equation}
 To estimate the other integral, suppose first that $1<p<2$. Using the fact that
 $X^2 \leq cX^{2-p}$ and (\ref{eq:1p2}) we have
 \begin{equation}
 \int_{\Omega}|\nabla v|^pd^{p-k}X^2dx \leq
 c\int_{\Omega}|\nabla v|^pd^{p-k}X^{2-p}dx
 \leq c I[u],
 \label{s3}
 \end{equation}
 and (\ref{eq:des}) follows from (\ref{s1}), (\ref{s2}), (\ref{s3}).

 Consider now the case $p \geq 2$. Using the fact that $X \leq 1$
 and  (\ref{eq:aa11}) we obtain in a similar fashion
 \begin{equation}
 \int_{\Omega}|\nabla v|^pd^{p-k}X^2dx \leq
 \int_{\Omega}|\nabla v|^pd^{p-k}Xdx
 \leq c I[u],
 \label{s4}
 \end{equation}
 and (\ref{eq:des}) follows from (\ref{s1}), (\ref{s2}), (\ref{s4}).$\hfill //$

 %
 %
 %
 %

\section{Appendix}
Here we present the two auxiliary lemmas that were used in Section 6. The
first one is a  one-dimensional Hardy type inequality.

 \begin{lemma}
 Let $p\in(1,\infty)$ and $q\geq p$ be given. For any $\alpha>-(p-1)$
 there exists $c>0$ such that
 \begin{equation}
 \int_0^1|v'|^pr^{p-1}X^{\alpha}dr \geq c
 \left(\int_0^1|v|^qr^{-1}X^{1+(\alpha +p-1)q/p}dr\right)^{p/q},
 \label{eq:th}
 \end{equation}
 for all $v\in C^{\infty}_c(0,1)$.
 \label{lem:4500}
 \end{lemma}

{\em Proof.} Apply [M, Theorem 3, p. 44] with 
$d\mu = r^{-1} X^{1+(\alpha +p-1)q/p} \chi_{[0,1]} dr$
and $d \nu = r^{p-1} X^{\alpha} \chi_{[0,1]} dr$. $\hfill //$

The second lemma is a weighted  Sobolev inequality.

\begin{lemma}
 Let $D>\sup_{x \in \Omega}d(x,K)$.
 Given $1<p<N$, $p<q\leq Np/(N-p)$ and $a\in\R$ there exists
 $c=C(p,q,D,\Omega)>0$ such that for all $v\in C^{\infty}_c(\Omega\setminus K)$
 \begin{eqnarray}
 &&\int_{\Omega}|v|^qd^{-N+(N-k)q/p}X^{\alpha q/p}(d/D)dx
 \leq \nonumber\\
 &\leq&c\left(
 \int_{\Omega}|\nabla v|^pd^{p-k}X^{\alpha}(d/D)dx +
 \int_{\Omega}|v|^pd^{-k}X^{\alpha}(d/D)dx
 \right)^{q/p}.
 \label{eq:rosa}
 \end{eqnarray}
 \label{lem:sobboun}
\end{lemma}
When $d$ is the distance from  the boundary $\partial \Omega$ (that is  $k=1$),
the above
result is given in [OK]; see Example 18.16  in p. 264 there. Since for the general
case we have not found a reference, we present a proof.

 {\em Proof. }Once again it suffices to consider the
 case $D=1$, the general case following by scaling.
 We shall make use of the standard Sobolev inequality
 \begin{equation}
 \int_{B(r)}|v|^qdx\leq cr^{N+q-Nq/p}\left(r^{-p}
 \int_{B(r)}|v|^pdx+\int_{B(r)}|\nabla v|^pdx\right)^{q/p},
 \quad v\in W^{1,p}(B(r)),
 \label{eq:zab}
 \end{equation}
 where $B(r)$ is any ball of radius $r$ and the constant is
 independent of $r$. Now, it follows from the
 Besicovich covering lemma (see \cite{M}) 
 that there exists a sequence $(x_m)$ of points
 in $\Omega$ with the following properties:
 defining, say,  $r_m=d(x_m)/10$, the balls
 \[B_m:=B(x_m,r_m),\quad m=1,2,\ldots\]
 satisfy
 \begin{eqnarray*}
 \ia && \Omega \subset \cup_m B_m;\\
 \ib && \mbox{there exists a number $M$ depending
 only on the dimension $N$}\\
 &&\mbox{such that each $x\in\Omega$
 belongs to at most $M$ of the $B_m$'s.}
 \end{eqnarray*}
 It follows from the choice of the radii $r_m$ that
 there exist constants $c_i,c_i'$ such that
 \begin{equation}
 \cases{c_1 r_m\leq d(x)\leq c_2 r_m, & \cr
 c_1' X(r_m)\leq X(d(x))\leq c_2' X(r_m), & }
 \;\; x\in B(x_m,r_m),\:\, m=1,2,\ldots .
 \label{eq:dr}
 \end{equation}
 This implies in particular that for any fixed $\theta,\eta\in\R$ we have
 \[
 c_1''r_m^{\theta}X^{\eta}(r_m)\int_{B_m}|u|^pdx\leq
 \int_{B_m}|u|^pd^{\theta}X^{\eta}(d)dx
 \leq c_2''r_m^{\theta}X^{\eta}(r_m)\int_{B_m}|u|^pdx,\]
 for all $m=1,2,\ldots$ and $u\in W^{1,p}(B_m)$.
 Hence
 \begin{eqnarray*}
 &&\hspace{-.8cm}
 \int_{\Omega}|v|^qd^{-N+(N-k)q/p}X^{\alpha q/p}(d)dx\\
 &\leq&\sum_{m=1}^{\infty}
 \int_{B_m}|v|^qd^{-N+(N-k)q/p}X^{\alpha q/p}(d)dx\\
 &\leq&c\sum_{m=1}^{\infty}r_m^{-N+(N-k)q/p}X^{\alpha q/p}(r_m)
 \int_{B_m}|v|^qdx\\
 &\leq&c\sum_{m=1}^{\infty}\left(\int_{B_m}|\nabla v|^pr_m^{p-k}
 X^{\alpha}(r_m)dx+\int_{B_m}|v|^pr_m^{-k}
 X^{\alpha}(r_m)dx\right)^{q/p}\\
 &\leq&c\sum_{m=1}^{\infty}\left(\int_{B_m}|\nabla v|^pd(x)^{p-k}
 X^{\alpha}(d(x))dx+\int_{B_m}|v|^pd(x)^{-k}
 X^{\alpha}(d(x))dx\right)^{q/p}\\
 &\leq&c\left(\sum_{m=1}^{\infty}\int_{B_m}|\nabla v|^pd(x)^{p-k}
 X^{\alpha}(d(x))dx+
 \sum_{m=1}^{\infty}\int_{B_m}|v|^pd(x)^{-k}
 X^{\alpha}(d(x))dx\right)^{q/p}\\
 &\leq&c\left(\int_{\Omega}|\nabla v|^pd(x)^{p-k}
 X^{\alpha}(d(x))dx+
 \int_{\Omega}|v|^pd(x)^{-k}
 X^{\alpha}(d(x))dx\right)^{q/p},
 \end{eqnarray*}
since by (ii) we have $\sum_{m}\int_{B_m}f\leq M\int_{\Omega}f$
for any non-negative function $f$ on $\Omega$.$\hfill //$





\begin{thebibliography}{RRR}

\bibitem[AS]{AS}{Ambrosio L. and Soner H.M. Level set approach to mean curvature flow in arbitrary
codimension. J. Diff. Geometry {\bf 43} (1996) 693-737.}

 \bibitem[BC]{BC}{Baras P. and Cohen L. Complete blow-up after $T_{{\rm max}}$ for the solution of
 a semilinear heat equation. J. Funct. Anal. {\bf 71} (1987) 142-174.}

\bibitem[BL]{BL}{Brezis H. and Lieb E.H. Sobolev inequalities with remainder terms. J. Funct. Anal.
{\bf 62} (1985) 73-86.}

 \bibitem[BM]{BM}{Brezis H. and Marcus M. Hardy's inequalities
 revisited. Ann. Scuola Norm. Pisa {\bf 25} (1997) 217-237.}

 \bibitem[BMS]{BMS}{Brezis H., Marcus M. and Shafrir I. Extremal functions for Hardy's inequality
 with weight. J. Funct. Anal. {\bf 171} (2000) 177-191.}


 \bibitem[BV]{BV}{Brezis H. and V\'{a}zquez J.-L. Blow-up solutions
 of some nonlinear elliptic problems. Rev. Mat. Univ. Comp. Madrid
 {\bf 10} (1997) 443-469.}


 \bibitem[CM]{CM}{Cabr\'{e} X. and Martel Y. Existence versus instantaneous blowup for linear heat
 equations with singular potentials. C.R. Acad. Sci. Paris Ser. I Math. {\bf 329} (1999) 973-978}

 \bibitem[D]{D}{Davies E.B. A review of Hardy inequalities. Oper. Theory Adv. Appl. {\bf 110}
(1998) 55-67.}

 \bibitem[DH]{DH}{Davies E.B. and Hinz A.M. Explicit constants for Rellich inequalities
 in $L_p(\Omega)$. Math. Z. {\bf 227} (1998) 511-523.}

 \bibitem[DM]{DM}{Davies E.B. and Mandouvalos N. The hyperbolic geometry and spectrum of irregular
domains. Nonlinearity {\bf 3} (1990) 913-945.}

 \bibitem[EG]{EG}{Evans L.C. and Gariepy R.F. Measure theory and fine properties
 of functions. CRC Press 1992.}

 \bibitem[FT]{FT}{Filippas S. and Tertikas A. Optimizing Improved
 Hardy inequalities. Preprint.}

\bibitem[FHT]{FHT}{Fleckinger J., Harrell II  M.  E. and  Thelin F.
Boundary behavior and estimates for solutions of equations containing the $p$-Laplacian.
Electron. J. Diff. Equns. {\bf 38} (1999) 1-19. }

\bibitem[GGM]{GGM}{Gazzola F., Grunau H.-Ch. and Mitidieri E. Hardy inequalities with optimal constants and
remainder terms. Preprint.}

 \bibitem[GP]{GP}{Garcia J.P. and Peral I. Hardy inequalities and some
 critical elliptic and parabolic problems. J. Diff. Equations {\bf 144}
 (1998) 441-476.}

 \bibitem[HKM]{HKM}{Heinonen J, Kilpelinen T. and Martio O. Nonlinear potential
theory of degenerate elliptic equations. The Clanderon Press, Oxford University Press,
New York, 1993.}

 \bibitem[HLP]{HLP}{Hardy G., P\'{o}lya G. and Littlewood J.E. Inequalities. 2nd edition,
 Cambridge University Press 1952.}

 \bibitem[L]{L}{Lindqvist P. On the equation $\diver(|\nabla u|^{p-2}
 \nabla u)+\lambda|u|^{p-2}u=0$. Proc. Amer. Math. Soc. {\bf 109} (1990) 157-164.}


 \bibitem[MMP]{MMP}{Marcus M., Mizel V.J. and Pinchover Y. On the best constant for Hardy's inequality
 in $\R^n$. Trans. Amer. Math. Soc. {\bf 350} (1998) 3237-3255.}


 \bibitem[MS]{MS}{Matskewich T. and Sobolevskii P.E. The best
 possible constant in generalized Hardy's inequality for
 convex domain in $\R^n$. Nonlinear Anal., Theory, Methods \& Appl.,
 {\bf 28} (1997) 1601-1610.}


 \bibitem[M]{M}{Maz'ja V.G. Sobolev spaces. Springer 1985.}

 \bibitem[OK]{OK}{Opic B. and Kufner A. Hardy-type inequalities. Pitman
 Research Notes in Math., vol.219, Longman 1990.}

 \bibitem[PV]{PV}{Peral I. and V\'{a}zquez J.L. On the stability or instability of the singular solution of the semilinear
 heat equation with exponential reaction term. Arch. Rational Mech. Anal. {\bf 129} (1995) 201-224.}

\bibitem[S]{S}{Sneider R. Convex bodies: The Brunn-Minkowski theory. Encyclopedia
of Math. and its Applications {\bf 44}, Cambridge University Press, 1993.}

 \bibitem[T]{T}{Tertikas A. Critical phenomena in linear elliptic problems. J. Funct. Anal. {\bf 154} (1998)
 42-66.}

 \bibitem[V]{V}{Vazquez J.L. Domain of existence and blowup for the exponential reaction-diffusion
 equation. Indiana Univ. Math. J. {\bf 48} (1999) 677-709.}



 \bibitem[VZ]{VZ}{Vazquez J.L. and Zuazua E., The Hardy inequality and
 the asymptotic behavior of the heat equation with an
 inverse-square potential. J. Funct. Anal. {\bf 173} (2000) 103-153.}

 \end{thebibliography}
 \end{document}